\documentclass{amsart}
\usepackage{amsmath,amsthm,amssymb,amsfonts,latexsym}
\usepackage{tikz}
\usepackage[colorlinks=true]{hyperref}

\hypersetup{urlcolor=blue, citecolor=blue}

\def\arXiv#1{   {\href{http://arxiv.org/abs/#1}
   {{\mdseries\ttfamily arXiv:#1}}}}

\def\doi#1{   {\href{http://dx.doi.org/#1}
   {{\mdseries\ttfamily DOI}}}}

\newcommand{\al}{\alpha}    \newcommand{\be}{\beta}
\newcommand{\de}{\delta}    
  \newcommand{\ep}{\varepsilon}
  \newcommand{\La}{\Lambda}

\newcommand{\ga}{\gamma}    \newcommand{\Ga}{\Gamma}
\newcommand{\R}{\mathbb{R}}

\newcommand{\pt}{\partial_t}\newcommand{\pa}{\partial}
\newcommand{\les}{{\lesssim}}
\newcommand{\beeq}{\begin{equation}}\newcommand{\eneq}{\end{equation}}

\newcommand{\Sp}{{\mathbb S}}

\newenvironment{prf}{\noindent {\bf Proof.} }{\endprf\par}
\def \endprf{\hfill  {\vrule height6pt width6pt depth0pt}\medskip}

\numberwithin{equation}{section}

\def\<{\langle}             \def\>{\rangle}	\def\({\left(}                 \def\){\right)}

\newtheorem{theorem}{Theorem}[section]
\newtheorem{lemma}[theorem]{Lemma}
\newtheorem{proposition}[theorem]{Proposition}

\theoremstyle{definition}

\theoremstyle{remark}
\newtheorem{remark}[theorem]{Remark}

\numberwithin{equation}{section}

\begin{document}
\title[Combined effects of two nonlinearities]
{Combined effects of two nonlinearities in lifespan of 
small solutions to semi-linear wave equations}

\author[K.\,Hidano]{Kunio Hidano}
\address{Department of Mathematics, Faculty of Education, 
Mie University, 1577 Kurima-machiya-cho Tsu, Mie Prefecture 514-8507, Japan}
\email{hidano@edu.mie-u.ac.jp}

\author[C.\,Wang]{Chengbo Wang}
\address{School of Mathematical Sciences,
                Zhejiang University,        Hangzhou 310027, China}
\email{wangcbo@zju.edu.cn}
\urladdr{http://www.math.zju.edu.cn/wang}

\author[K.\,Yokoyama]{Kazuyoshi Yokoyama}
\address{Hokkaido University of Science, 
7-15-4-1 Maeda, Teine-ku, Sapporo, Hokkaido 006-8585, Japan}
\email{yokoyama@hus.ac.jp}

\thanks{
The authors are very grateful to the referees for their helpful comments.
The first author was supported in part by 
the Grant-in-Aid for Scientific Research (C) (No.\,23540198), 
Japan Society for the Promotion of Science (JSPS). 
The second author was supported by 
NSFC 11301478, 11271322, and
Zhejiang Provincial Natural Science Foundation of 
China LR12A01002.}



\date{\today}

\dedicatory{In memory of Rentaro Agemi}



\begin{abstract}
This paper investigates the combined effects of two 
distinctive power-type nonlinear terms (with parameters $p,q>1$) in  the 
lifespan of small solutions to semi-linear wave equations. 
We determine the full region of $(p,q)$ to admit global existence of small solutions, at least for spatial dimensions $n=2, 3$.  Moreover, for many $(p,q)$ when there is no global existence, we obtain sharp lower bound of the lifespan, which is of the same order as the upper bound of the lifespan.
\end{abstract}

\maketitle

\section{Introduction}

In this paper, we are interested in determining the dichotomy between the global solvability and the blow up, for a large class of the small-amplitude semilinear wave equations with two distinctive power-type nonlinear terms.
More precisely, letting $p, q>1$, and considering the sample wave equations with parameters $(p,q)$ 
\begin{equation}\label{nlw}
\partial_t^2 u-\Delta u=|\partial_t u|^p+|u|^q,\,\,\,
t>0,\,\,x\in{\mathbb R}^n,
\end{equation}
we are interested in determining the region of $(p,q)$,  for which
the following statement is true:
for any given nontrivial pair of compactly supported smooth functions $(f,g)$, there exists a small parameter $\ep_0=\ep_0(f,g,n,p,q)>0$, such that the problem with initial data of sufficiently small size $\ep\in (0,\ep_0)$
\begin{equation}\label{eq-data0}
u(0,x)=\varepsilon f(x),\,\,\,
\partial_t u(0,x)=\varepsilon g(x),
\end{equation}
 admits global solutions.
For the cases where there is no global existence, 
we are also interested in the estimate of lifespan, denoted by $T_\ep(p,q)$, from above and below, in terms of the parameters
$(\ep, n, p,q)$.  Here, the lifespan is defined as the supremum of $T>0$ such that the problem admits a unique solution in $[0,T]\times\R^{n}$.

When the spatial dimension is one, the standard ordinary differential inequality argument (see, e.g., \cite{Kato80,Si83,Si84,YZ,ZhouHan}) could be easily adapted to show that, for any $p, q>1$, the problem does not admit global solutions in general. On the basis of this fact, in what follows, we will assume $n\ge 2$. 

The problem can be regarded as a natural combination of the following two well-investigated problems
\begin{equation}\label{nlw1S}
\partial_t^2 v-\Delta v=|v|^q,\,\,\,t>0,\,\,x\in{\mathbb R}^n\ ,
\end{equation}
\begin{equation}\label{nlw2G}
\partial_t^2 w-\Delta w=|\pt w|^p,\,\,\,t>0,\,\,x\in{\mathbb R}^n\ .
\end{equation}
The first problem \eqref{nlw1S} is in relation with the Strauss conjecture, for which the critical
power, denoted by $q_c(n)$, is known to be
the positive root of the quadratic equation \begin{equation}
(n-1)q^2-(n+1)q-2=0\eneq
that is,
\beeq
q_c(n)
:=
\frac{n+1+\sqrt{n^2+10n-7}}{2(n-1)}\ .
\end{equation}
This problem was initiated in \cite{Jo}, where the critical value was determined to be $p_c=1+\sqrt{2}$ for $n=3$. 
Shortly afterward, \cite{Str81} included the conjecture that the critical power is given by $q_c(n)$.
The existence portion of the conjecture was verified in
\cite{Gl2} ($n=2$), \cite{Zh1995} ($n=4$), \cite{LS} ($n\le 8$),  \cite{KK1995, LS, KK1998} ($n\geq 5$ and radial data) and  \cite{GLS, Ta} (all $n\ge 2$).
 The necessity of $q > q_c$ for small data global existence is due to
\cite{Jo, Gl1,Si84,Sc,YZ, Zh07}.
Concerning the estimates of the lifespan, denoted by $T_\ep^{S}$, for $1<q\le q_c$,
it is known that, for some constant $C>0$, we have
(see \cite{TaWa11, ZhHa14-bu} and references therein)
$$T_\ep^{S}\le S_\ep(q):=\left\{\begin{array}{ll}
 C \ep^{-\frac{2q(q-1)}{2(q+1)-(n-1)q(q-1)}},& q<q_c,\\
 \exp(C \ep^{-q(q-1)}),& q=q_c,
 \end{array}\right.$$
 for any $\ep\in (0,1)$,
 which is known to be sharp at least for $\max(1, 2/(n-1))<q<q_c$, or $q=q_c$ with $n\le 8$ (see \cite{LaZ, LS, Ta15} and references therein). See \cite{HMSSZ, SSW, LMSTW, Wa15p} and references therein for recent works for this problem on various spacetimes.
 
 Concerning \eqref{nlw2G}, it is conjectured that the critical power is given by 
 $$(n-1)(p_c-1)=2,\ p_c(n)=1+\frac{2}{n-1}\ ,$$
 which is known as the Glassey conjecture.
For the case $1<p\le p_c$, nonexistence of global small solutions and upper bound of the lifespan (denoted by $T_{\ep}^{G}$)
 $$T_\ep^{G}\le G_\ep(p):=\left\{\begin{array}{ll}
 C \ep^{-\frac{2(p-1)}{2-(n-1)(p-1)}},& p<p_c,\\
 \exp(C \ep^{-(p-1)}),& p=p_c,
 \end{array}\right.$$
has been known through the works
\cite{Jo81, Sc86, Ra, Ag, Zh}.
For the existence part with $p>p_c$, global existence of small solutions has been proved for $n=2, 3$ in \cite{HT,Tz} (with an earlier radial $3$-D result \cite{Si83}). 
Recently, the present authors have succeeded in extending the results 
of \cite{Si83,HT,Tz} to the case of higher space dimensions 
$n\geq 4$ under the radial assumption of the initial data \cite{HWY}. The existence part of the Glassey conjecture with general data for $n\ge 4$ remains unsolved.
See \cite{Wa14, Wa15} for recent works for this problem on various spacetimes.
 
\begin{figure}
\centering
\begin{tikzpicture} [scale=1.5]
\filldraw[black!20!white] (1,1)--(1,3.97)--(2,3.97)--(2,3) to[out=297,in=134] (2.414,2.414)--(5,2.414)--(5,1)--(1,1);
\draw[thick, dashed, domain=2:2.414] plot (\x, {{1+2/((\x)-1)}}) ;
\node[right] at (2.2,2.8) {$(q-1)(p-\frac{2}{n-1})=\frac{4}{n-1}$};
\node[below] at (2.414,2.414) {$(q_c,q_c)$}; 
\draw[fill=black,line width=1pt] (2,3) circle[radius=0.3mm];
\node[left] at (1,3) {$1+\frac{4}{n-1}$};
\draw[fill=black,line width=1pt] (2.414,2.414) circle[radius=0.3mm];
\node[below left] at (1,1) {$(1,1)$};
\node[left] at (1,2.414) {$q_c$};
\node[below] at (2.414,1) {$q_c$};
\node[below] at (2,1) {$p_c$};
	\draw[thick] (2,3)--(2,3.97) (2.414,2.414)--(5,2.414);
	\draw[thick,-stealth] (1,1)--(1,4) node[left]{$q$};
	\draw[thick,-stealth] (1,1)--(5,1) node[below]{$p$};
\node at (3.3,3.3) {Global?};
\end{tikzpicture}
\caption{Blow up region for $(p,q)$: the shaded region (except the broken curve)}
\label{3D-bu}
\end{figure}
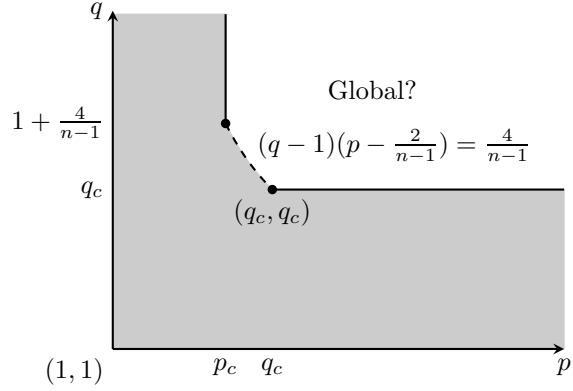

On the basis of the known blow up results for the Strauss conjecture and the Glassey conjecture, 
by comparing the nonlinearities, it would not be difficult to adapt the proof to conclude that
\beeq\label{region-SG}
 T_\ep(p,q)\le \left\{\begin{array}{ll}
 S_\ep(q), &q\le q_c\ ,\\
 G_\ep(p),& p\le p_c\ .\end{array}\right.\eneq  (For the definition of $T_\varepsilon(p,q)$, see the first paragraph of this section.) Then to admit global small solutions, we are forced to consider $q>q_c$ and $p>p_c$.
Recently, Han and Zhou \cite{HZ} studied the blow up  phenomenon for \eqref{nlw}. Among others, they found a new  combined effect 
on the lifespan, by proving the blow up results for $n\ge 2$, $q>q_c$, $p>p_c$ and 
\begin{equation}\label{region-ZH}
(q-1)((n-1)p-2)<4.
\end{equation} Moreover, they  obtained an upper bound of the lifespan
\begin{equation}\label{life-ZH}
T_\varepsilon(p,q)
\leq Z_\ep(p,q):=
C\varepsilon^{-\frac{2p(q-1)}{2(q+1)-(n-1)p(q-1)}},
\end{equation} 
for the cases where $$\max(1,\frac{2}{n-1})<p\le \frac{2n}{n-1},\ 1<q<\frac{2n}{n-2},\ (q-1)((n-1)p-2)<4.$$
Here, by checking the proof, we observe that, the restriction $p>2/(n-1)$ is not necessary. That is,
the upper bound \eqref{life-ZH} is
actually valid for
$$1<p\le \frac{2n}{n-1},\ 1<q<\frac{2n}{n-2},\ (q-1)((n-1)p-2)<4.$$
 See Figure \ref{3D-bu} for the illustration of the known region of non-existence of global small solutions. It is interesting to observe that the conformal power $1+4/(n-1)$ occurs very naturally in the figure.
 
 It will then be very natural and interesting to ask whether \eqref{region-SG} and \eqref{region-ZH} are the only regions of blow up or not, and
try to determine the full region of global existence. 
Of course, we expect global existence when $p,q$ are both large enough (recall that we are considering small data problem). In the case where one of the powers (say $q$) is relatively small comparing the other, the nonlinearity with that power will tend to be dominant and one may infer the  behavior of problem similar to the problem with only one power-type nonlinear term (say $|u_t|^p+|u|^q\sim |u|^q$ for the case $q\ll p$). The essential difficulty comes from the case where neither of the powers is large.

Assuming the blow up results \eqref{region-SG}, \eqref{region-ZH} have been precise enough, with the experiences from the Strauss conjecture and the Glassey conjecture, we may naturally infer that we have global existence for 
\begin{equation}\label{region-glo}
q>q_c, p>p_c,
(q-1)((n-1)p-2)>4\ ,
\eneq
and may only admit almost global existence for the ``critical" case
\begin{equation}\label{region-glo2}
q>q_c, p>p_c,
(q-1)((n-1)p-2)=4\ .
\eneq
Surprisingly enough, in this paper, we are able to prove global existence, not only for \eqref{region-glo}, but also for the ``critical" case \eqref{region-glo2}. This is our first main theorem.
As the Glassey conjecture for $n\ge 4$, with general data, remains open, it is very natural for us to restrict ourselves to the case
 $n=2, 3$ in our main theorems.

\begin{theorem}\label{thm-1}
Let $n=2, 3$,
$s_d:=1/2-1/q$,
\begin{equation}
q>q_c,\,\,\,p>p_c\,\,\mbox{and}\,\,\,
(q-1)((n-1)p-2)\ge 4.
\end{equation}
Suppose that $f\in \dot H^1\cap \dot H^{s_d}$ and $g\in L^2\cap \dot H^{s_d-1}$ with \beeq\label{eq-data1}\La:=\sum_{|\be|\le 2, |\alpha|\le \min(2,|\be|+1)}
\bigl(
\|x^\al\nabla_x^\be f\|_{{\dot H}^1\cap{\dot H}^{s_d}}
+
\|x^\al\nabla_x^\be g\|_{L^2\cap{\dot H}^{s_d-1}}\bigl)<\infty.\eneq
Then, there exists an $\varepsilon_0>0$ depending on 
$n$, $p$, $q$, and $\Lambda$ such that 
the Cauchy problem \eqref{nlw}-\eqref{eq-data0} admits a unique global solution,
provided that $\varepsilon\in [0,\varepsilon_0)$. 
\end{theorem}

For a more precise statement of the result, see Theorem \ref{thm-1-2}.
See Figure \ref{3D-bu2} for the region of global existence, that is, the shaded region (except the broken lines).
Here, and in what follows, 
by 
$\|\cdots\|_{{\dot H}^1\cap{\dot H}^{s_d}}$
we naturally mean
\begin{equation}
\|\varphi\|_{{\dot H}^1\cap{\dot H}^{s_d}}
:=
\|\varphi\|_{{\dot H}^1({\mathbb R}^n)}
+
\|\varphi\|_{{\dot H}^{s_d}({\mathbb R}^n)}.
\end{equation}

 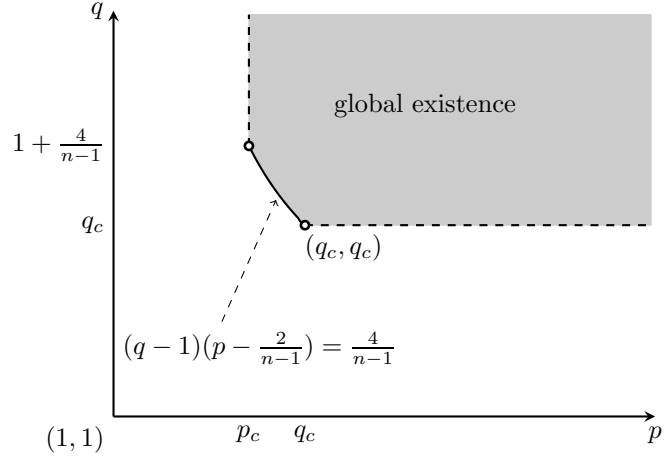
\begin{figure}
\centering
\begin{tikzpicture} [scale=1.8]
\filldraw[black!20!white] (4.97,3.97)--(2,3.97)--(2,3) to[out=297,in=134] (2.414,2.414)--(4.97,2.414)--(4.97,3.97);
\draw[thick, domain=2:2.414] plot (\x, {{1+2/((\x)-1)}}) ;
\node[below] at (2.7,2.414) {$(q_c,q_c)$}; 
\node[left] at (1,3) {$1+\frac{4}{n-1}$};

\node[below left] at (1,1) {$(1,1)$};
\node[left] at (1,2.414) {$q_c$};
\node[below] at (2.414,1) {$q_c$};
\node[below] at (2,1) {$p_c$};
	\draw[thick, dashed] (2,3)--(2,3.97) (2.414,2.414)--(4.97,2.414);
\draw[fill=white,line width=1pt] (2.414,2.414) circle[radius=0.3mm];
\draw[fill=white,line width=1pt] (2,3) circle[radius=0.3mm];
	\draw[thick,-stealth] (1,1)--(1,4) node[left]{$q$};
	\draw[thick,-stealth] (1,1)--(5,1) node[below]{$p$};
\node at (3.3,3.3) {global existence};
\draw[->, dashed] (1.8,1.7)--(2.2,2.6);
\node[right] at (1,1.5) {$(q-1)(p-\frac{2}{n-1})=\frac{4}{n-1}$};

\end{tikzpicture}
\caption{Region to admit global existence: the shaded region (except the broken lines)}
\label{3D-bu2}
\end{figure}

For the cases where there is no global existence, that is, when
$q\le q_c$, or $p\le p_c$, or $(q-1)((n-1)p-2)<4$,
we are also interested in the estimate of lifespan $T_\ep$, from above and below. 
Observe that
for the cases where $1<p\le 2n/(n-1)$, $1<q< q_c$ and $(q-1)((n-1)p-2)<4$, we have
$$S_\ep \ge Z_\ep \ \mathrm{ for }\ \ep\ll 1 \ \Leftrightarrow q\ge p\ ,$$
and
for the cases where $1<p< p_c$, $1<q<2n/(n-2)$ and $(q-1)((n-1)p-2)<4$,
we have
$$G_\ep \ge Z_\ep \ \mathrm{ for }\ \ep\ll 1 \Leftrightarrow q\le 2p-1\ .$$
Then we know  
from 
\eqref{region-SG} and
\eqref{life-ZH} that
$$T_\ep\le \left\{
\begin{array}{ll}
G_\ep (p),     &    1< p\le p_c,\ q\ge 2p-1,\\
S_\ep (q),     &    1<q\le q_c,\ q\le p,\\
Z_\ep (p,q),  & (q-1)((n-1)p-2)<4,\ 2p-1 \ge q\ge p>1
\end{array}\right.
$$
when $n\ge 2$.
On the basis of these observations, it is natural to infer that the sharp lower bound is of the same size as these upper bounds.
In the following second main theorem, we obtain the sharp lifespan estimates, for
$n=2, 3$,  $q>2/(n-1)$ and $q, p\ge 2$,
except the critical case $q=q_c\le p$.
\begin{theorem}\label{thm-2}
Let $n=2, 3$, 
 $q>2/(n-1)$ and $q, p\ge 2$. Assume also $q\le q_c$, $p\le p_c$ or $(q-1)((n-1)p-2)<4$.
Then for any
$f\in \dot H^1\cap \dot H^{s_d}$ and $g\in L^2\cap \dot H^{s_d-1}$ with $\La<\infty$,
there exists an $\varepsilon_0>0$ depending on 
$n$, $p$, $q$, and $\Lambda$ such that 
the Cauchy problem \eqref{nlw}-\eqref{eq-data0} admits a unique solution for $t\in [0, T]$,
provided that $\varepsilon\in (0,\varepsilon_0)$, where
$$T= \left\{
\begin{array}{ll}
G_\ep (p),     &    2\le p\le p_c,\ q\ge 2p-1,\\
S_\ep (q) ,    &    2\le q< q_c,\ 2/(n-1)<q\le p,\\
\exp(c\ep^{1-q}),     &    q=q_c\le p,\\
Z_\ep (p,q),  & (q-1)((n-1)p-2)<4,\ 2\le p\le q\le 2p-1
\end{array}\right.
$$
for some small constant $c>0$.
Moreover, the lower bound of the lifespan is sharp in general, except $q=q_c\le p$.
\end{theorem}
For a more precise statement of the result, see Theorem \ref{thm-2-2}.
See Figures \ref{3D} and \ref{2D} for the illustration of the lower bound of the lifespan. 
The existence result for smooth nonlinearity with $(p,q,n)=(3,4,2)$ and $T=Z_{\ep}(3,4)$ was known from \cite{Ka01}.
 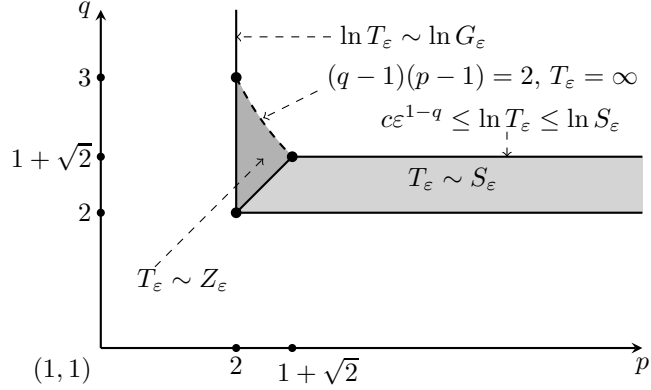
\begin{figure}
\centering
\begin{tikzpicture} [scale=1.8]

\filldraw[black!17!white] (2,2)--(2.414,2.414)--(5,2.414)--(5,2)--(2,2);
\filldraw[black!30!white] (2,2)--(2,3) to[out=297,in=134] (2.414,2.414)--(2,2);

\draw[thick, dashed, domain=2:2.414] plot (\x, {{1+2/((\x)-1)}}) ;
\draw[thick] (2,2)--(2.414,2.414) (2,2)--(2,3);
\draw[thick] (2,2)--(5,2);


\draw[fill=black,line width=1pt] (2,3) circle[radius=0.3mm];
\draw[fill=black,line width=1pt] (2.414,2.414) circle[radius=0.3mm];

\node[below left] at (1,1) {$(1,1)$};

\node[left] at (1,3) {$3$};
\draw[fill=black,line width=1pt] (1,3) circle[radius=0.2mm];
\node[left] at (1,2.414) {$1+\sqrt{2}$};
\draw[fill=black,line width=1pt] (1,2.414) circle[radius=0.2mm];
\node[below] at (2.6,1) {$1+\sqrt{2}$};
\draw[fill=black,line width=1pt] (2.414, 1) circle[radius=0.2mm];
\node[below] at (2,1) {$2$};
\draw[fill=black,line width=1pt] (2, 1) circle[radius=0.2mm];
\node[left] at (1,2) {$2$};
\draw[fill=black,line width=1pt] (1, 2) circle[radius=0.2mm];

	\draw[thick] (2,3)--(2,3.5) (2.414,2.414)--(5,2.414);
	\draw[thick,-stealth] (1,1)--(1,3.5) node[left]{$q$};
	\draw[thick,-stealth] (1,1)--(5,1) node[below]{$p$};


\node[right] at (2.7,3.3) {$\ln T_\ep\sim \ln G_\ep$};
	\draw[->, dashed] (2.7,3.3)--(2, 3.3);

\node[right] at (3, 2.7) {$c \ep^{1-q}\le \ln T_\ep\le \ln S_\ep$};
	\draw[<-, dashed] (4,2.414)--(4, 2.6);

\node[right] at (3.2, 2.25) {$T_\ep\sim  S_\ep$};

\node[right] at (1.2,1.5) {$T_\ep\sim Z_\ep$};
	\draw[->, dashed] (1.4,1.6)--(2.2, 2.4);
	
	\node[right] at (2.6,3) {$(q-1)(p-1)=2$, $T_\ep=\infty$};
	\draw[->, dashed] (2.6,2.9)--(2.2, 2.7);

\draw[fill=black,line width=1pt] (2,2) circle[radius=0.3mm];
\end{tikzpicture}
\caption{Dimension three: estimates of the lifespan 
}
\label{3D}
\end{figure}

 \begin{figure}
\centering
\begin{tikzpicture} [scale=1.5]

\filldraw[black!17!white] (2,2)--(3.56,3.56)--(6,3.56)--(6,2)--(2,2)
(2,3)--(3,5)--(3,6)--(2,6)--(2,3);
\filldraw[black!25!white] (2,2)--(2,3)--(3,5)
 to[out=285,in=120](3.56,3.56)--(2,2);

\draw[thick, dashed, domain=3:3.56] plot (\x, {{1+4/((\x)-2)}}) ;
\draw[thick, dashed] (2,2)--(6,2);


\draw[fill=black,line width=1pt] (3,5) circle[radius=0.3mm];
\draw[fill=black,line width=1pt] (3.56,3.56) circle[radius=0.3mm];
\draw[fill=black,line width=1pt] (2, 3) circle[radius=0.3mm];

\node[below left] at (1,1) {$(1,1)$};

\node[left] at (1,5) {$5$};
\draw[fill=black,line width=1pt] (1,5) circle[radius=0.2mm];
\node[left] at (1,3) {$3$};
\draw[fill=black,line width=1pt] (1,3) circle[radius=0.2mm];
\node[left] at (1,3.56) {$(3+\sqrt{17})/2$};
\draw[fill=black,line width=1pt] (1,3.56) circle[radius=0.2mm];

\node[below] at (3.76, 1) {$(3+\sqrt{17})/2$};
\draw[fill=black,line width=1pt] (3.56, 1) circle[radius=0.2mm];
\node[below] at (3,1) {$3$};
\draw[fill=black,line width=1pt] (3, 1) circle[radius=0.2mm];
\node[left] at (1,3) {$3$};
\draw[fill=black,line width=1pt] (1, 3) circle[radius=0.2mm];
\node[below] at (2,1) {$2$};
\draw[fill=black,line width=1pt] (2, 1) circle[radius=0.2mm];
\node[left] at (1,2) {$2$};
\draw[fill=black,line width=1pt] (1, 2) circle[radius=0.2mm];

	\draw[thick] (2,2)--(2,3)--(3,5)--(3,6) (2,6)--(2,3) (2,2)--(3.56,3.56)--(6,3.56);
	\draw[thick,-stealth] (1,1)--(1,6) node[left]{$q$};
	\draw[thick,-stealth] (1,1)--(6,1) node[below]{$p$};


\node[right] at (3.8,5.5) {$\ln T_\ep\sim \ln G_\ep$};
	\draw[->, dashed] (3.7,5.5)--(3, 5.5);

\node[right] at (3.8, 4) {$c \ep^{1-q}\le \ln T_\ep\le \ln S_\ep$};
	\draw[<-, dashed] (5,3.56)--(5, 3.9);

\node[right] at (4, 3) {$T_\ep\sim  S_\ep$};

\node[right] at (2,5.3) {$T_\ep\sim G_\ep$};

\node[right] at (2.3,3.5) {$T_\ep\sim Z_\ep$};
	
	\node[right] at (3.35,4.4) {$(q-1)(p-2)=4$, $T_\ep=\infty$};
	\draw[->, dashed] (4,4.3)--(3.4, 4);

	\node at (3,1.65) {$q=p$};
	\draw[->, dashed] (3,1.7)--(3, 2.95);

	\node[right] at (3.5,4.75) {$q=2p-1$};
	\draw[->, dashed] (3.5,4.65)--(2.5, 3.9);

\draw[fill=white,line width=1pt] (2,2) circle[radius=0.3mm];
\end{tikzpicture}
\caption{Dimension two: estimates of the lifespan (except the broken line)}
\label{2D}
\end{figure}

\begin{remark}
For the critical case $q=q_c\le p$ and $n=2,3$, we know that
$$\exp(c\ep^{1-q})\le T_\ep\le \exp(C \ep^{-q(q-1)})\ .$$
It will be interesting to determine the sharp estimate of the lifespan, and we infer that the sharp estimate will be $T_\ep\ge  \exp(c \ep^{-q(q-1)})$ for $q=q_c\le p$.
\end{remark}

\begin{remark}
As is obvious from the proofs, 
Theorems \ref{thm-1} and \ref{thm-2} remain valid 
for the equations of the form
$$
\partial_t^2 u
-
\Delta u
=
F_{q}(u)
+
G_{p}(\pt u,\nabla u),
$$
where $F_{q}$ and $G_{p}$ are $C^{2}$ functions with 
$$|\pa_{u}^{j}F_{q}(u)|\le C_{j} |u|^{q-j}, \forall u\in [-1,1], j=0,1,2,$$
$$\sum_{|\al|=j}|\pa_{v}^{\al}G_{p}(v)|\le C_{j} |v|^{p-j}, \forall v\in \R^{n+1}, |v|\le 1, |\al|\le 2.$$
\end{remark}
Let us conclude the introduction by describing the strategy of the proof.
It is natural to view the problem  \eqref{nlw} as either a perturbation of \eqref{nlw1S}
by a forcing term $|\pt u|^p$, or a perturbation of \eqref{nlw2G}
by a forcing term $|u|^q$. One of the remarkable difficulties of the problem lies in the fact that the standard proofs for these two problems are typically distinct, which forces us to seek a robust proof of existence results which could handle both of the nonlinearities effectively.
Fortunately, there does exist a method of proof which has the same nature and works well for both of the nonlinearities. 
 
 Actually,  inspired by \cite{LZ}, the first author \cite{Hi} developed an alternative proof of global existence of small solutions to \eqref{nlw1S} with
$q>q_c(n)$ and $n=2, 3, 4$, by using the homogeneous Sobolev space ${\dot H}^{s_d}$ in the iteration argument. It turns out that such a method works well also for the \eqref{nlw2G} when $n=2, 3$, by using the standard energy space $\dot H^1\times L^2$. As the proof for \eqref{nlw1S} is more involved than that for \eqref{nlw2G}, to prove Theorem \ref{thm-1} and Theorem \ref{thm-2},
it is natural to view the problem as a perturbation of \eqref{nlw1S} by a forcing term $|\pt u|^p$.
 

The equation \eqref{nlw} has the ``forcing term'' 
$|\partial _t u|^p$, 
which involves a higher-order derivative of $u$. 
This naturally leads us to a modification of the norm 
in the iteration scheme, 
and we allow for some growth of the ${\dot H}^{s_d}$ norm,
see \eqref{grow} and \eqref{grow2}.  

This paper is organized as follows. 
In the next section, we collect several preliminary inequalities.
In Section \ref{sec-3}, we give the setup for the existence results, and prove the basic iteration estimates, which is the key in the proof of Theorems \ref{thm-1} and \ref{thm-2}. 
Then, using the key estimates obtained in Section \ref{sec-3}, 
we present  the proof of global existence, Theorem \ref{thm-1}, and
long time existence, Theorem \ref{thm-2},
in Sections \ref{sec-4} and \ref{sec-5}, on a case-by-case basis.

\subsection{Notation}
For $x=(x_1, \dots, x_n)$,  we will use
polar coordinates $x=r\omega$ with $r=|x|$, $\omega\in\mathbb{S}^{n-1}$,
and the full space-time gradient $\partial=(\pt,\nabla_x)=(\partial_0,\partial_1,\dots ,\partial_n)$.
In addition, we denote
$L_j=t\partial_j+x_j\partial_0$ ($1\le j\le n$), 
$\Omega_{kl}=x_k\partial_l-x_l\partial_k$ $(1\leq k<l\leq n)$, 
$L_0=t\partial_0+x \cdot \nabla$. 
The collection of all these operators is denoted by $\Gamma_j$, with
$0\le j\le \nu:=(n^2+3n+2)/2$. 
For a multi-index $\alpha=(\alpha_0,\dots,\alpha_\nu)$, 
$\Gamma^\alpha:=\Gamma^{\alpha_0}_0\cdots\Gamma^{\alpha_\nu}_\nu$. 
Moreover we will employ the notation 
$\langle x\rangle:=\sqrt{1+|x|^2}$ for $x\in \R^n$,
$\|\Ga^{\le k} u\|:=\sum_{|\al|\le k}\|\Ga^{\al} u\|$,
and use the Fourier multiplier $|D|:=\sqrt{-\Delta}$. 

The homogeneous space $\dot H^s$ for $|s|<n/2$ is the completion of Schwartz functions with respect to the norm $\||D|^s u\|_{L^2}$.
We will use the following mixed-norm $L^{q_1}_t L^{q_2}_r L^{q_3}_\omega$,
$$\|f\|_{L^{q_1}_t L^{q_2}_r L^{q_3}_\omega}=\left\|\left(\int_{0}^\infty \|f(t, r\omega)\|_{L^{q_3}_\omega}^{q_2} r^{n-1} dr\right)^{1/q_2}\right\|_{L^{q_1}(t> 0)},$$ with trivial modification for the case $q_2=\infty$,
where $L^q_\omega$ is the standard Lebesgue space on the sphere $\Sp^{n-1}$. Occasionally, we will omit the subscripts. Also, at times we will employ abbreviations, such as  $L^{q_1} L^{q_2}=L^{q_1}_t  L^{q_2}_r L^{q_2}_\omega$ and $L^q_T=L^q([0,T])$.

Let $2\le p<\infty$, $2\le q<\infty$, $q> 2/(n-1)$, $s_d:=1/2- 1/q \in [0,1/2)\cap ((2-n)/2, 1/2)$.
We introduce
$$X_u^k(t):= \|\Ga^{\le k}u(t)\|_{\dot H^{s_d}},\ Y_u^k(t):=\|\pa \Ga^{\le k} u(t)\|_{L^2},\ 
Z_{u}^k(t):=X_u^k(t)^{\frac{q+1}{q+2}}Y_{u}^k(t)^{\frac 1{q+2}}\ ,
$$
with abbreviation $X_u=X_u^0$, $Y_u=Y_u^0$, $Z_u=Z_u^0$,
and define $p_1, p_2\in (1,2), q_1, q_2\in [2,\infty]$ as follows
\begin{equation}\label{eq-def-p_1}
\frac n{p_1}=1-s_d+\frac n2,\ 
\frac{n-1}{p_2}=
\frac{n}{2}-s_d=
\frac n{q_1}, \ \frac{1}{q_1}+\frac{1}{q_2}=\frac12,
\end{equation}
such that we have the following
$$L^{p_1}\subset \dot H^{s_d-1},\ L^{p_2}_\omega\subset H^{s_d-1/2}_\omega,\ 
\dot H^{s_d}\subset L^{q_1 }.
$$
Observe also that we have
\beeq\label{eq-rel}
\frac{1}{p_1}=\frac{1}{n}+\frac{1}{q_1}=\frac{1}{n}+\frac{1}{2}-\frac 1{q_2}, \ \frac{1}{p_2}=\frac{1}{q(n-1)}+\frac{1}{2}, q_2=\frac{n}{s_d}\ .\eneq

Let $\chi_1(t,x)$ be the characteristic function of the set 
$\{x\in{\mathbb R}^n:|x|<(1+|t|)/2\}$ and $\chi_2=1-\chi_1$. 
We will also use $A\les B$ to stand for $A\le C B$ where the constant $C$ may change from line to line.
In addition,
when denoting by $a+$ (or $a-$) for $a\in{\mathbb R}$, 
we mean that the relevant estimate holds for $a+\varepsilon$ (or $a-\varepsilon$) for sufficiently small $\varepsilon>0$. Also, the notation $\infty-$ means that the relevant estimate holds for sufficiently large values.

\section{Preliminaries}
In this section, we give some preliminary results.

\begin{proposition}[Sobolev inequalities]\label{thm-Sobo}
For any $s\in [0, n/2)$, we have
\begin{equation}\label{2.1}
\|v\|_{L^{q_0}({\mathbb R}^n)}\les \|v\|_{{\dot H}^s({\mathbb R}^n)},  \|v\|_{{\dot H}^{-s}({\mathbb R}^n)}\les \|v\|_{L^{q_0'}({\mathbb R}^n)},\ \frac{n}{q_0}=\frac{n}2-s.
\end{equation}
\end{proposition}

\begin{proposition}\label{thm-trace-crit} Let $n\ge 2$, $s_d=1/2-1/q$ with $2\le q<\infty$. Then 
we have
\begin{equation}\label{2.2}
\|r^{(n-1)s_d}v\|_{L^q_{r} L^2_\omega}\les \|v\|_{{\dot H}^{s_d}}\ .\end{equation}
\end{proposition}

\begin{prf}
See Theorem 2.10 of \cite{LZ}. 
For the reader's convenience, 
we give an alternative proof of \eqref{2.2},
by using the endpoint trace inequality and the real interpolation. 
Recall the following  endpoint trace inequality (see \cite{FW})
\begin{equation}
\sup_{r>0}
r^{(n-1)/2}\|u(r\cdot)\|_{L^2(\Sp^{n-1})}
\les \|u\|_{{\dot B}^{1/2}_{2,1}}, n\ge 2\ .
\end{equation}
Here, and in the following discussion, 
by ${\dot B}^s_{p,q}={\dot B}^s_{p,q}({\mathbb R}^n)$ 
we mean the homogeneous Besov space, see, e.g., Chapter 6 of \cite{BL}. 
Observe also the obvious equality
\begin{equation}
\|r^{(n-1)/2}u(r\cdot)\|_{L^2({\mathbb R}^+;L^2(\Sp^{n-1}))}
=
\|u\|_{L^2({\mathbb R}^n)}\ ,
\end{equation}
where ${\mathbb R}^+:=(0,\infty)$. 
Let $T(u):=r^{(n-1)/2}\|u(r\cdot)\|_{L^2(\Sp^{n-1})}$.
Then,
by Theorems 3.1.2, 6.4.5, and 5.2.1 of \cite{BL} 
(see also {\it Remark} on page 41 of \cite{BL}), 
we see that the sublinear operator 
$$
T:{\dot B}^{1/2}_{2,1}({\mathbb R}^n)+L^2({\mathbb R}^n)
\to L^\infty({\mathbb R}^+)+L^2({\mathbb R}^+)
$$
satisfies 
$$
T:[{\dot B}^{1/2}_{2,1}({\mathbb R}^n),\,{\dot B}^0_{2,2}({\mathbb R}^n)]
_{2/q,q}
\to
[L^\infty({\mathbb R}^+),\,L^2({\mathbb R}^+)]_{2/q,q}
$$
and so is the inequality 
\begin{equation}
\|T(u)\|_{L^q({\mathbb R}^+)}
\les \|u\|_{{\dot B}^{1/2-1/q}_{2,q}({\mathbb R}^n)}
\les\|u\|_{{\dot H}^{1/2-1/q}({\mathbb R}^n)}
, 2<q<\infty\ 
\end{equation}
which gives us \eqref{2.2}.
\end{prf}

\begin{proposition}\label{thm-Klai-Sobo}
If $1\leq p<\infty$ and $s>n/p$, then 
the inequality
\begin{equation}\label{2.3}
(1+|t|+|x|)^{(n-1)/p}(1+||t|-|x||)^{1/p}|v(t,x)|
\les \|\Ga^{\le s} v(t,\cdot)\|_{L^p_x}
\end{equation}
holds. 
If $1\leq p<q<\infty$ and $1/q\geq 1/p-1/n$, 
then we have
\begin{equation}\label{2.4}
\|\chi_1 v(t,\cdot)\|_{L^q_x}
\les (1+|t|)^{-n(1/p-1/q)}
\|\Ga^{\le 1}v(t,\cdot)\|_{L^p_x}\ .
\end{equation}
\end{proposition}
 See \cite{Kl87} and \cite[Theorem 2.9]{LZ} for
the proof of \eqref{2.3} and
 \eqref{2.4}.
 
We will also need the following trace estimates.
\begin{proposition}\label{thm-trace}
Let $n\geq 2$. 
Then 
the inequalities 
\begin{equation}\label{2.5}
\|r^{(n/2)-s}
v(r\omega)\|_{L^\infty_r H^{s-1/2}_\omega}
\les \|v\|_{{\dot H}^s({\mathbb R}^n)}, \frac 12<s<\frac n2
\end{equation}
\begin{equation}\label{2.6}
\|r^{(n/2)-s}
v(r\omega)\|_{L^\infty_r L^{p}_\omega}
\les
\|v\|_{{\dot H}^s({\mathbb R}^n)}, \frac 12<s<\frac n2, \frac{n-1}p=\frac n2-s.
\end{equation}
hold.
In addition, if $2\le p\leq 4$, $q=2p/(4-p)$,
then we have
\begin{equation}\label{eq-trace-ed}
r^{(n-1)/2}\|u(r\cdot)\|_{L^p_\omega}
\leq
\sqrt{p}
\|\partial_r u\|_{L^2}^{1/2}
\|u\|_{L^2_rL^q_\omega}^{1/2}.
\end{equation}
If 
$p\in [2, \min\left\{4, 2(n-1)/(n-2)\right\} ]$ except the endpoint $p=4$ and $n=3$,
then 
\begin{equation}\label{2.9}
\sup_{r>0}r^{(n-1)/2}
\|v(r\cdot)\|_{L^p_\omega}
\leq
C
\|\partial_r v\|_{L^2({\mathbb R}^n)}^{1/2}
\|\Omega^{\le 1} v\|_{L^2({\mathbb R}^n)}^{1/2}\ .
\end{equation}
\end{proposition}

\begin{prf} For the proof of the trace lemma \eqref{2.5}, 
see \cite{Ho} for $n\geq 3$ and \cite{FW} for $n\geq 2$. 
By the Sobolev embedding on the unit sphere $\Sp^{n-1}$, 
we obtain \eqref{2.6} directly from \eqref{2.5}. 

Let us turn to the inequality \eqref{eq-trace-ed},
which generalizes the well-known inequality
\begin{equation}\label{eq-tr-prf}
r^{(n-1)/2}\|u(r\cdot)\|_{L^2_\omega}
\leq
\sqrt2
\|\partial_r u\|_{L^2({\mathbb R}^n)}^{1/2}
\|u\|_{L^2({\mathbb R}^n)}^{1/2}\ .
\eneq
It suffices to give the proof for $u\in C_0^\infty({\mathbb R}^n)$. 
We use a natural modification of the proof of \eqref{eq-tr-prf}. 
We first note for any fixed $R>0$
\begin{eqnarray}\label{eq-tr-prf1}
& &
\bigl(
R^{(n-1)/2}\|u(R\cdot)\|_{L^p_\omega}
\bigr)^p
=
\int_{\Sp^{n-1}}
R^{(n-1)p/2}
|u(R\omega)|^pd\omega\\
& &
\leq
p\int_R^\infty\int_{\Sp^{n-1}}
r^{(n-1)p/2}|u(r\omega)|^{p-1}|(\omega\cdot\nabla u)(r\omega)|drd\omega
\nonumber\\
& &
\leq
p\|r^{(n-1)\theta/2}u\|_{L^{2(p-1)}({\mathbb R}^n)}^{p-1}
\|\partial_r u\|_{L^2({\mathbb R}^n)},\nonumber
\end{eqnarray}
where $\theta:=(p-2)/(p-1)$. 
Using $1/(2(p-1))=\theta/p+(1-\theta)/q$,
$(1-\theta)(p-1)=1$, we get
\beeq\label{eq-tr-prf2}
\|r^{(n-1)\theta/2}u\|_{L^{2(p-1)}({\mathbb R}^n)}
\leq
\|
r^{(n-1)/2}u\|_{L^\infty_r L^p_\omega}
^\theta
\|u\|_{L^{2}_r L^q_\omega}^{1-\theta}.\eneq
The inequalities \eqref{eq-tr-prf1}-\eqref{eq-tr-prf2}  yield
$$
\|
r^{(n-1)/2}u\|_{L^\infty_r L^p_\omega}^p
\leq
p
\|
r^{(n-1)/2}u\|_{L^\infty_r L^p_\omega}^{p-2}
\|u\|_{L^2_rL^q_\omega}
\|\partial_r u\|_{L^2},$$
which gives \eqref{eq-trace-ed}.

Then \eqref{2.9}
follows  immediately from \eqref{eq-trace-ed}, if we recall
the embedding $H^1_\omega \hookrightarrow L^q_\omega$ 
with $q=2p/(4-p)$,
for $p\in [2, \min\left\{4, 2(n-1)/(n-2)\right\} ]$ except the endpoint case $p=4$ and $n=3$.
\end{prf}

By Sobolev inequality \eqref{2.1} and duality to the trace lemma \eqref{2.6}, we have for $2/(n-1)<q<\infty$, i.e., $1-s_d\in (1/2,n/2)$,
\beeq\label{eq-SoboTr}\|F\|_{\dot H^{s_d-1}}
\les \|\chi_1 F\|_{L^{p_1}}+
\| \<t\>^{-(n-2)/2-s_d}
\chi_2 F\|_{ L^1_r L^{p_2}_\omega }
\ ,\eneq
where $p_1$ and $p_2$ are defined in
\eqref{eq-def-p_1}.
Thus, using the standard energy estimates, we get the following
\begin{proposition}\label{thm-ener}
Let $n\geq 2$ and $2/(n-1)<q<\infty$. 
Then, for any $T>0$ we have
\beeq\label{2.10}
\|\pa u(T)\|_{{\dot H}^{s_d-1}}
\les
\|\pa u(0)\|_{{\dot H}^{s_d-1}}
+
 \|\chi_1 F\|_{L^1_T  L^{p_1}}+
\|
\<t\>^{-(n-1)/2+1/q}
\chi_2 F\|_{L^1_T L^1_r L^{p_2}_\omega}
\eneq
and
\begin{equation}\label{2.10-ener}
\|\pa u(T)\|_{L^2}
\les
\|u(0)\|_{{\dot H}^1}+\|\pt u(0)\|_{L^2}
+
\|F\|_{L^1_T L^2},
\end{equation}
for any solutions to the inhomogeneous wave equation
$\partial_t^2 u-\Delta u=F$.
\end{proposition}

%
Moreover, we have the following classical relations by direct computations. 
\begin{proposition}\label{thm-commu}
The following commuting relations hold$:$
\begin{eqnarray}
& &
[\Gamma_i,\Box]=0\,\,\,\mbox{for $i=0,\dots,\nu-1$, and}\,\,\,
[L_0,\Box]=-2\Box,\\
&&[\pa_j,\Gamma_k]=\sum_{l=0}^n C^{j,k}_l\pa_l,
[\Gamma_j,\Gamma_k]=\sum_{l=0}^\nu C^{j,k}_l\Gamma_l,\,\,\,
j,\,k=0,\dots,\nu.
\end{eqnarray}
Here $C^{j,k}_l$ denotes a constant depending on 
$j$, $k$, and $l$.
\end{proposition}


In particular, we see by this proposition the equivalence 
between $\|\Ga^{\le 2}\partial u(t,\cdot)\|_{L^2_x}$ and 
$\|\partial \Gamma^{\le 2} u(t,\cdot)\|_{L^2_x}$. 
This fact will be repeatedly employed. 

\section{Iteration}\label{sec-3}
Let $n= 2, 3$, $q\in [2,\infty)\cap (2/(n-1),\infty)$, $p\in [2,\infty)$,
$s_d=1/2-1/q$,
 $(f,g)$ be functions with
\eqref{eq-data1}
and the solution space \beeq\label{3.2ST}
S_T= \{u: 
\pa \Gamma^\alpha u\in C([0,T];{\dot H}^{s_d-1}\cap L^2),\,
|\alpha|\leq 2,
u(0)=\ep f, \pt u(0)=\ep g
\}\eneq
with $S_\infty=\cap_{T>0} S_T$.
For any $u\in S_T$,
we define $P u$ as the solution of the following linear wave equation
\beeq\label{nlw-Pi}(\pt^2-\Delta) P u=F(u):=a |u|^q+b |\partial_t u|^p, Pu(0)=\ep f, \pt (Pu)(0)=\ep g\ ,\eneq and then solving \eqref{nlw} in $S_T$ is equivalent to finding the fixed point $u$ such that $Pu=u$.

In view of Propositions \ref{thm-ener} and \ref{thm-commu}, to show $Pu\in S_T$ for $u\in S_T$, we need to obtain the initial bound on 
$(\Gamma^\alpha u,\partial_t\Gamma^\alpha u)_{|t=0}$, which is ensured by the following proposition.
\begin{proposition}[Initial data]\label{thm-data}
Let $n=2, 3$, $u\in S_T$.
Rewrite $(\Gamma^\alpha P u)(0)$ and 
$(\partial_t\Gamma^\alpha P u)(0)$ in terms of 
$f$ and $g$, through the equation \eqref{nlw-Pi}.
Then there exists a constant $C_0>0$ depending only on $n$, $p$, $q$, and 
independent of $\varepsilon\in (0,1)$, such that
\beeq\label{eq-data}
\sum_{|\alpha|\leq 2}
\bigl(
\|\pa \Gamma^\alpha P u(0)\|_{{\dot H}^{s_d-1}}
+
\|\pa \Gamma^\alpha P u(0)\|_{L^2}\bigr)
\leq
C_0M\varepsilon\ .\eneq
 Here, we have defined 
\begin{equation}
M:=\Lambda+\Lambda^{2p-1}+\Lambda^{p+q-1}.
\end{equation}
\end{proposition}

As is standard, to prove the existence of the fixed point, we typically show that the operator $P$ maps certain ball in $S_T$ into itself, and satisfies certain (weak) contraction property on that ball.
\begin{proposition}[Uniform boundedness]\label{thm-unifbd}
Let $n= 2, 3$, $q\in [2,\infty)\cap (2/(n-1),\infty)$, $p\in [2,\infty)$, $s=s_d$.
There exists a constant $C_1\ge C_0$ such that we have for any $u\in S_T$ with $T\in (0,\infty)$,
\beeq\label{eq-key-1}
X_{Pu}^2(T)
\le C_1 M\ep+C_1\int_0^T \left(\<t\>^{-(q-1)(\frac n2-s)+1}X_u^2(t)^q + \<t\>^{-(p-1)\frac{n-1}2+\frac 12-s}Y_u^2(t)^p \right) dt\ ,\eneq
\begin{eqnarray}
Y_{Pu}^2(T)& 
\le  & C_1 M\ep+C_1 \int_0^T \left(\<t\>^{-(p-1)\frac{n-1}2}Y_u^2(t)^p + \<t\>^{-q(\frac n2-s)+\frac n 2}
X_u^2(t)^q\right)dt \nonumber\\
&&+C_1\int_0^T \<t\>^{-(q-1)\frac{n-1}2}Z_u^2(t)^qdt\ .
\label{eq-key-2}
\end{eqnarray}
\end{proposition}

\begin{remark}Setting
$$\al_1:=-(q-1)(\frac n2-s)+1=\frac 1q-\frac{n-1}2(q-1),$$
$$\al_2:=-(p-1)\frac{n-1}2+\frac 12-s=\frac 1q-\frac{n-1}2(p-1),$$
$$\al_3:=-(p-1)\frac{n-1}2,\ 
\al_4:=-(q-1)\frac{n-1}2,$$
$$\al_5:=-q(\frac n2-s)+\frac n 2=-\frac 12-\frac{n-1}2(q-1),$$
we can rephrase \eqref{eq-key-1}-\eqref{eq-key-2} as follows
\beeq\label{eq-key-1'}
X_{Pu}^2(T)
\le C_1 M\ep+C_1\int_0^T \left(\<t\>^{\al_1}X_{u}^2(t)^q + \<t\>^{\al_2}Y_{u}^2(t)^p \right) dt\ ,\eneq
\beeq
Y_{Pu}^2(T)
\le  C_1 M\ep+C_1 \int_0^T \left( \<t\>^{\al_3}Y_{u}^2(t)^p +\<t\>^{\al_4}X_u^2(t)^{\frac{q+1}{q+2}q}Y_u^2(t)^{\frac{q}{q+2}}+ \<t\>^{\al_5}
X_u^2(t)^q 
\right)
dt\ ,
\label{eq-key-2'}
\eneq
for any $T>0$.
\end{remark}

\begin{proposition}[Convergence]\label{thm-conv}
Let $n= 2, 3$, $q\in [2,\infty)\cap (2/(n-1),\infty)$, $p\in [2,\infty)$, $s=s_d$. 
Then there exists a constant $C_2\ge C_1$ such that for any $T\in (0,\infty)$, we have for any $u, v\in S_T$,
\begin{eqnarray}
X_{Pu-Pv}(T)&
\le&C_2 \int_0^T \<t\>^{\al_1}(X_u^2(t)+X_v^2(t))^{q-1}X_{u-v}(t) dt\nonumber\\
&& + C_2\int_0^T \<t\>^{\al_2}(Y_u^2(t)+Y_v^2(t))^{p-1} Y_{u-v}(t)dt\ ,
\label{eq-key-3}
\end{eqnarray}
\begin{eqnarray}
Y_{Pu-Pv}(T)& 
\le &C_2 \int_0^T \<t\>^{\al_3}(Y_u^2(t)+Y_v^2(t))^{p-1}Y_{u-v}(t) dt \nonumber\\
&&
+C_2\int_0^T \<t\>^{\al_4}(Z_u^2(t)+Z_v^2(t))^{q-1}Z_{u-v}(t)dt \nonumber\\
&&
+C_2\int_0^T \<t\>^{\al_5}
(X_u^2(t)+X_v^2(t))^{q-1}X_{u-v}(t) dt
\ .
\label{eq-key-4}
\end{eqnarray}
\end{proposition}

\subsection{Proof of Propositions \ref{thm-unifbd} and \ref{thm-conv}}
In this subsection, we present the proof of Propositions \ref{thm-unifbd} and \ref{thm-conv}, which could be viewed as the main estimates of this paper. 
We will put off the tedious proof of Proposition \ref{thm-data} till the next subsection.

By Propositions \ref{thm-ener}, \ref{thm-commu}  and \ref{thm-data}, for the proof of  Proposition \ref{thm-unifbd},
we need only to prove that for fixed $t>0$,
\beeq\label{eq-claim1}
 \|\chi_1 \Ga^{\le 2} |u(t,\cdot)|^q\|_{L^{p_1}}+
\|
\<t\>^{-(n-1)/2+1/q}
\chi_2  \Ga^{\le 2} |u|^q\|_{L^1_r L^{p_2}_\omega}
\les \<t\>^{\al_1}X_{u}^2(t)^q ,\eneq
\beeq\label{eq-claim2}
 \|\chi_1 \Ga^{\le 2} |\pt u|^p\|_{L^{p_1}}+
\|
\<t\>^{-(n-1)/2+1/q}
\chi_2  \Ga^{\le 2} |\pt u|^p\|_{L^1_r L^{p_2}_\omega}
\les \<t\>^{\al_2}Y_{u}^2(t)^p ,\eneq
\beeq\label{eq-claim3}
 \|\Ga^{\le 2} |\pt u|^p\|_{L^{2}}
\les \<t\>^{\al_3}Y_{u}^2(t)^p ,\eneq
\beeq\label{eq-claim4}
 \|\Ga^{\le 2} | u|^q\|_{L^{2}}
\les \<t\>^{\al_4}X_u^2(t)^{\frac{q+1}{q+2}q}Y_u^2(t)^{\frac{q}{q+2}}+ \<t\>^{\al_5}
X_u^2(t)^q\ .\eneq
Similarly, 
 the proof of Proposition \ref{thm-conv} is reduced to
\beeq\label{eq-claim1'}
\||u|^q-|v|^q\|_{\dot H^{s_d-1}}
\les \<t\>^{\al_1}(X_{u}^2(t)+X_v^2(t))^{q-1}X_{u-v}(t),\eneq
\beeq\label{eq-claim2'}
\||\pt u|^p-|\pt v|^p\|_{\dot H^{s_d-1}}
\les \<t\>^{\al_2}(Y_{u}^2(t)+Y_{v}^2(t))^{p-1}Y_{u-v}(t),\eneq
\beeq\label{eq-claim3'}
 \| |\pt u|^p-|\pt v|^p\|_{L^{2}}
\les \<t\>^{\al_3}(Y_{u}^2(t)+Y_{v}^2(t))^{p-1}Y_{u-v}(t),\eneq
\beeq\label{eq-claim4'}
 \| | u|^q-|v|^q\|_{L^{2}}
\les \<t\>^{\al_4}(Z_u^2(t)+Z_u^2(t))^{q-1}Z_{u-v}(t)+ \<t\>^{\al_5}
(X_{u}^2(t)+X_v^2(t))^{q-1}X_{u-v}(t)\ .\eneq

\subsubsection{Estimates \eqref{eq-claim1} and \eqref{eq-claim1'}}
We start with the estimate of \eqref{eq-claim1}, for which
we have $$|\Ga^{\le 2} |u|^q|\les |u|^{q-2} |\Gamma^{\le 1} u|^2+
 |u|^{q-1} |\Gamma^{\le 2} u|\ .$$

Using \eqref{2.3}-\eqref{2.4} and \eqref{2.1}, we get
\begin{eqnarray}
\|\chi_1 |u|^{q-2} |\Gamma^{\le 1} u|^2\|_{L^{p_1}}
&\leq &
\|\chi_1 u\|_{L^\infty}^{q-2}
\|\chi_1 \Ga^{\le 1} u\|_{L^{2p_1}}^2\nonumber\\
& \les&
\langle t\rangle^{-\frac{n}{q_1}(q-2)-n(\frac{1}{q_1}-\frac{1}{2p_1})\times2}
\| \Ga^{\le 2} u\|_{L^{q_1}}^q
\les
\langle t\rangle^{\al_1}X_u^2(t)^q.
\end{eqnarray}
In the case of $2p_1<q_1$, we have actually employed 
the H\"older inequality to get 
$\|\chi_1 v\|_{L^{2p_1}}
\leq
C\<t\>^{-n(1/q_1-1/(2p_1))}
\|\chi_1 v\|_{L^{q_1}}$. 
Using $1/p_1=1/n+1/q_1$, we also get
\begin{eqnarray}
\|\chi_1 |u|^{q-1}\Gamma^{\le 2} u\|_{L^{p_1}}
&\leq&
\|\chi_1 u\|_{L^{n(q-1)}}^{q-1}
\|\chi_1 \Ga^{\le 2} u\|_{L^{q_1}}\nonumber\\
& \les&
\<t\>^{-n(\frac{1}{q_1}-\frac 1{n(q-1)})(q-1)}
\|\Ga^{\le 2} u\|_{L^{q_1}}^q
\les
\< t\>^{\al_1}X_u^2(t)^q.
\end{eqnarray}
The same kind of proof gives us
\begin{eqnarray}
\|\chi_1 (|u|^{q}-|v|^q)\|_{L^{p_1}}
&\les&
(\|\chi_1 u\|_{L^{n(q-1)}}+\|\chi_1 v\|_{L^{n(q-1)}})^{q-1}
\|\chi_1 (u-v)\|_{L^{q_1}}\nonumber\\
& \les&
\< t\>^{\al_1}(X_u^2(t)+X_v^2(t))^{q-1}X_{u-v}(t).\label{eq-claim1'-prf}
\end{eqnarray}

Using 
the Sobolev embedding on $\Sp^{n-1}$, $H^1_\omega\subset L^{\infty-}_\omega$, $H^2_\omega\subset L^\infty_\omega$, and \eqref{2.2}, 
we obtain
\begin{eqnarray}
& &
\|\chi_2 |u|^{q-2}|\Gamma^{\le 1} u|^2\|_{L^1_r L^{p_2}_\omega}
\leq
\|\chi_2 u\|_{L^q_rL^\infty_\omega}^{q-2}
\|\chi_2  \Ga^{\le 1} u\|_{L^{q}_rL^{2p_2}_\omega}^2\\
& &
\les
\|\chi_2 \Ga^{\le 2} u\|_{L^{q}_rL^{2}_\omega}^q
\les
\<t\>^{-(n-1)(\frac 12-\frac 1q)q}X_u^2(t)^q=\< t\>^{\al_1-(1/q-(n-1)/2)}X_u^2(t)^q.\nonumber
\end{eqnarray}
As $p_2< 2$,
we also obtain
\begin{eqnarray}
\|\chi_2 |u|^{q-1}\Gamma^{\le 2} u\|_{L^1_r L^{p_2}_\omega}
&\les& 
\|\chi_2 u\|_{L^q_r L^{\infty-}_\omega}^{q-1}
\|\chi_2 \Ga^{\le 2} u\|_{
L^{q}_rL^{2}_\omega
}\les
\|\chi_2 \Ga^{\le 1} u\|_{
L^{q}_rL^{2}_\omega}^{q-1}
\|\chi_2 \Ga^{\le 2} u\|_{
L^{q}_rL^{2}_\omega}
\nonumber\\
&
\les &
\<t\>^{-(n-1)(\frac 12-\frac 1q)q}X_u^2(t)^q=\< t\>^{\al_1-(1/q-(n-1)/2)}X_u^2(t)^q,
\end{eqnarray}
which completes the proof of \eqref{eq-claim1}.

Similarly, we have
\begin{eqnarray}
\|\chi_2 (|u|^q-|v|^q)\|_{L^1_r L^{p_2}_\omega}
&\les& 
(\|\chi_2 u\|_{L^q_r L^{\infty-}_\omega}+\|\chi_2 v\|_{L^q_r L^{\infty-}_\omega})^{q-1}
\|\chi_2 (u-v)\|_{
L^{q}_rL^{2}_\omega
}
\nonumber\\
&
\les &
\< t\>^{\al_1-(1/q-(n-1)/2)}(X_u^2(t)+X_v^2(t))^{q-1}X_{u-v}(t),
\end{eqnarray}
which, together with \eqref{eq-claim1'-prf}, gives  \eqref{eq-claim1'}.

\subsubsection{Estimates \eqref{eq-claim2} and \eqref{eq-claim2'}}
Next, we deal with \eqref{eq-claim2},
 for which
we have $|\Ga^{\le 2} |\pt u|^p|\les |\pt u|^{p-2} |\Gamma^{\le 1} \pt u|^2+
 |\pt u|^{p-1} |\Gamma^{\le 2} \pt u|$.
Using \eqref{2.3}-\eqref{2.4}, 
we get
\begin{eqnarray}
\|\chi_1 |\pt u|^{p-2} |\Gamma^{\le 1} \pt u|^2\|_{L^{p_1}}
&\leq &
\|\chi_1 \pt u\|_{L^\infty}^{p-2}
\|\chi_1 \Ga^{\le 1}\pt  u\|_{L^{2p_1}}^2\nonumber\\
& \les&
\langle t\rangle^{-\frac{n}{2}(p-2)-n(\frac{1}{2}-\frac{1}{2p_1})\times2}
\| \Ga^{\le 2}\pt  u\|_{L^{2}}^p\nonumber\\
&=&
\langle t\rangle^{-n(p-1)/2+1-s_d}
\| \Ga^{\le 2}\pt  u\|_{L^{2}}^p\nonumber\\
&\les&
\langle t\rangle^{\al_2}Y_u^2(t)^p,
\end{eqnarray}
where we have used the assumption $p\ge 2$ in the last inequality.
Similarly, using $1/p_1=(q+2)/(2nq)+1/2$, we also get
\begin{eqnarray}
\|\chi_1 |
\pt u|^{p-1}\Gamma^{\le 2} \pt  u\|_{L^{p_1}}
&\leq&
\|\chi_1 \pt u\|_{L^{\infty}}^{p-2}
\|\chi_1 \pt u\|_{L^{\frac{2nq}{q+2}}}
\|\chi_1 \Ga^{\le 2}\pt  u\|_{L^{2}}\nonumber\\
& \les&
\<t\>^{-(n/2)(p-2)-n(1/2-(q+2)/(2nq))}
\|\Ga^{\le 2} \pt u\|_{L^{2}}^p\nonumber\\
&=&
\langle t\rangle^{-n(p-1)/2+1-s_d}
\| \Ga^{\le 2}\pt  u\|_{L^{2}}^p\nonumber\\
&\les&
\langle t\rangle^{\al_2}Y_u^2(t)^p.
\end{eqnarray}
In the same vein, we could also get
\begin{eqnarray}
&&\|\chi_1 (|\pt u|^{p}-|\pt v|^{p})\|_{L^{p_1}}\nonumber\\
&\les&
(\|\chi_1 \pt u\|_{L^{\infty}}^{p-2}
\|\chi_1 \pt u\|_{L^{\frac{2nq}{q+2}}}+
\|\chi_1 \pt v\|_{L^{\infty}}^{p-2}
\|\chi_1 \pt v\|_{L^{\frac{2nq}{q+2}}})
\|\chi_1 \pt  (u-v)\|_{L^{2}}\nonumber\\
& \les&
\langle t\rangle^{\al_2}(Y_u^2(t)+Y_v^2(t))^{p-1}Y_{u-v}(t).\label{eq-claim2'-prf}
\end{eqnarray}

On the other hand, employing the Sobolev embedding on $\Sp^{n-1}$ 
and \eqref{2.9}, we get
\begin{eqnarray}
& &
\|\chi_2 |\pt u|^{p-2} |\Gamma^{\le 1} \pt u|^2\|
_{L^1_r L^{p_2}_\omega}
\leq
\|\chi_2 \partial_t u\|_{L^\infty}^{p-2}
\|\chi_2 \Ga^{\le 1}\partial_t u\|_{L^2_r L^{2 p_2}_\omega}^2\\
& &
\les
\|\chi_2 
\Ga^{\le 1}\partial_t u\|_{L^\infty_r L^{2+}_\omega}^{p-2}
\| \Ga^{\le 2}\partial_t u\|_{L^2}^2
\les
\<t\>^{-(n-1)(p-2)/2}\|\Ga^{\le 2}\partial_t u\|_{L^2}^p.\nonumber
\end{eqnarray}
Using $1/p_2=1/2+1/(q(n-1))$, 
we also obtain
\begin{eqnarray}
& &
\|\chi_2 |
\pt u|^{p-1}\Gamma^{\le 2} \pt  u\|_{L^1_r L^{p_2}_\omega}
\leq
\|\chi_2 \partial_t u\|_{L^\infty}^{p-2}
\| \partial_t u\|_{L^2_r L^{q(n-1)}_\omega}
\|\Ga^{\le 2}\partial_t u\|_{L^2}
\\
& &
\les
\<t\>^{-(n-1)(p-2)/2}
\|\Ga^{\le 2}\partial_t u\|_{L^2}^p.\nonumber
\end{eqnarray}
Observing that
$-(n-1)(p-2)/2=\al_2-(1/q-(n-1)/2)$, we see that this completes the proof of \eqref{eq-claim2}.

Similarly,
\begin{eqnarray}
& &
\|\chi_2  (|\pt u|^{p}-|\pt v|^{p})\|_{L^1_r L^{p_2}_\omega}\\
& \les&
(\|\chi_2 \partial_t u\|_{L^\infty}^{p-2}
\| \partial_t u\|_{L^2_r L^{q(n-1)}_\omega}
+\|\chi_2 \partial_t v\|_{L^\infty}^{p-2}
\| \partial_t v\|_{L^2_r L^{q(n-1)}_\omega})
\|\partial_t (u-v)\|_{L^2}
\nonumber\\
&\les &
\<t\>^{\al_2-(1/q-(n-1)/2)}(Y_u^2(t)+Y_v^2(t))^{p-1}Y_{u-v}(t)
,\nonumber
\end{eqnarray}
which, combined with \eqref{eq-claim2'-prf}, gives \eqref{eq-claim2'}.

\subsubsection{Estimates \eqref{eq-claim3} and \eqref{eq-claim3'}}


For the estimate of $|\pt u|^{p-2} |\Gamma^{\le 1} \pt u|^2$, 
we get by using \eqref{2.3}-\eqref{2.4}
$$
\|\chi_1 |\pt u|^{p-2} |\Gamma^{\le 1} \pt u|^2\|
_{L^2}\le
\|\chi_1\partial_t u\|_{L^\infty}^{p-2}
\|\chi_1\Ga^{\le 1}\partial_t u\|_{L^4}^2
\les
\langle t\rangle^{-(n/2)(p-1)}
\|\Ga^{\le 2}\partial_t u\|_{L^2}^p\ ,$$
and using \eqref{2.3}, \eqref{2.9}, and the Sobolev embedding 
on $\Sp^{n-1}$, 
\begin{eqnarray}
& &
\|\chi_2 |\pt u|^{p-2} |\Gamma^{\le 1} \pt u|^2\|_{L^2}\\
& &
\leq
\|\partial_t u\|_{L^\infty}^{p-2}
\|\chi_2\Ga^{\le 1}\partial_t u\|_{L^\infty_r L^{2+}_\omega}
\|\Ga^{\le 1} \partial_t u\|_{L^2_r L^{\infty-}_\omega}
\les
\langle t\rangle^{-(n-1)(p-1)/2}
\|\Ga^{\le 2}\partial_t u\|_{L^2}^p\ .\nonumber
\end{eqnarray}
In addition, by \eqref{2.3}, we get
\beeq
\||\partial_t u|^{p-1}\Gamma^{\le 2}\partial_t u\|_{L^2}
\le
\|\partial_t u\|_{L^\infty}^{p-1}
\|\Ga^{\le 2}\partial_t u\|_{L^2}
\les
\langle t\rangle^{-(n-1)(p-1)/2}
\|\Ga^{\le 2}\partial_t u\|_{L^2}^p.\eneq

Similarly, we obtain
\begin{eqnarray}
\| |\pt u(t, \cdot)|^p-|\pt v(t, \cdot)|^p\|_{L^{2}}
 & \les &  (\| \pt u\|_{L^\infty}+\| \pt v\|_{L^\infty})^{p-1}\|\pt (u-v)\|_{L^{2}}\nonumber\\
 & \les & \langle t\rangle^{-(n-1)(p-1)/2}
(Y_u^2(t)+Y_v^2(t))^{p-1}\|\pt (u-v)\|_{L^{2}}\ .
\label{eq-4.1'}
\end{eqnarray}

\subsubsection{Estimates \eqref{eq-claim4} and \eqref{eq-claim4'}}

For the estimate of $\|\chi_1 \Ga^{\le 2}|u(t,\cdot)|^q\|_{L^2}$, 
we get by using \eqref{2.3}-\eqref{2.4} and \eqref{2.1}
\begin{eqnarray}
\|\chi_1 |u|^{q-2}|\Gamma^{\le 1} u|^2\|_{L^2}
&\leq&
\|\chi_1 u\|_{L^\infty}^{q-2}
\|\chi_1 \Ga^{\le 1} u\|_{L^4}^2 \nonumber\\
& \les&
\large(
\langle t\rangle^{-n/q_1}
\|\Ga^{\le 2}u\|_{L^{q_1}}
\large)^{q-2}
\large(
\langle t\rangle^{-n(1/q_1-1/4)}
\|\Ga^{\le 2} u\|_{L^{q_1}}
\large)^2\nonumber\\
&\les &
\langle t\rangle^{-(n/q_1)q+(n/2)}X_u^2(t)^q
=\langle t\rangle^{\al_5}X_u^2(t)^q
\end{eqnarray}
and
\begin{eqnarray}
\|\chi_1 |u|^{q-1}\Gamma^{\le 2} u\|_{L^2}
&\leq
 &
\|\chi_1 u\|_{L^{q_2(q-1)}}^{q-1}
\|\Ga^{\le 2}u\|_{L^{q_1}}\nonumber\\
& \les&
\langle t\rangle^{-n(\frac{1}{q_1}-\frac{1}{q_2(q-1)})(q-1)}
\|\Ga^{\le 2}u\|_{L^{q_1}}^q\nonumber\\
&\les&
\langle t\rangle^{\al_5}X_u^2(t)^q.
\end{eqnarray}
Similarly, we have
\begin{eqnarray}
\|\chi_1 (|u|^{q}-|v|^q)\|_{L^2}
&\les
 &
( \|\chi_1 u\|_{L^{q_2(q-1)}}^{q-1}
+ \|\chi_1 v\|_{L^{q_2(q-1)}}^{q-1}
)
\|u-v\|_{L^{q_1}}\nonumber\\
&\les&
\langle t\rangle^{\al_5}(X_u^2(t)+X_v^2(t))^{q-1}X_{u-v}(t).\label{eq-claim4'-prf}
\end{eqnarray}

For the estimate of $\|\chi_2 \Ga^{\le 2}|u(t,\cdot)|^q\|_{L^2}$, 
we proceed as follows. 
Noting $1/2=(q-2)/(2q)+1/q$ 
and 
$1/2-1/(2q)=s_d(1-\theta)+\theta$ 
for $\theta:=1/(q+2)$, 
we obtain by using the Sobolev embedding on $\Sp^{n-1}$ 
and \eqref{2.2} 
\begin{eqnarray}
\|\chi_2 |u|^{q-2}|\Gamma^{\le 1} u|^2\|_{L^2}& \les &
\|\chi_2 u\|_{L^{2q}_rL^\infty_\omega}^{q-2}
\|\chi_2 \Ga^{\le 1}u\|_{L^{2q}_r L^4_\omega}^2\nonumber\\
&\les&
\|\chi_2\Ga^{\le 2}u\|_{L^{2q}_rL^2_\omega}^q\nonumber\\
& 
\les
&
\langle t\rangle^{-(n-1)(1/2-1/(2q))q}
\||D|^{1/2-1/(2q)}\Gamma^{\le 2} u\|_{L^2}^q\nonumber\\
& 
\les&
\langle t\rangle^{-(n-1)(q-1)/2}
\|\Gamma^{\le 2} u\|_{\dot H^{s_d}}^{q(1-\theta)}
\|\Gamma^{\le 2} u\|_{\dot H^1}^{q \theta }
\nonumber\\
& 
\les
&
\langle t\rangle^{-(n-1)(q-1)/2}
X_u^2(t)^{q\frac{q+1}{q+2}}
Y_u^2(t)^{\frac{q}{q+2}}=\<t\>^{\al_4}Z_u^2(t)^q
.\nonumber
\end{eqnarray}
Moreover, 
noting $1/2=(q-1)/(2q)+1/(2q)$ 
and using the Sobolev embedding on $\Sp^{n-1}$ and \eqref{2.2}, we obtain
\begin{eqnarray}
\|\chi_2 |u|^{q-1} \Gamma^{\le 2} u\|_{L^2}& \les &
\|\chi_2 u\|_{L^{2q}_rL^\infty_\omega}^{q-1}
\|\chi_2 \Ga^{\le 2}u\|_{L^{2q}_r L^2_\omega}\nonumber\\
&\les&
\|\chi_2\Ga^{\le 2}u\|_{L^{2q}_rL^2_\omega}^q
\les \<t\>^{\al_4}Z_u^2(t)^q
,\nonumber
\end{eqnarray}
and this completes the proof of \eqref{eq-claim4}.
Similarly, we have
\begin{eqnarray}
\|\chi_2 (|u|^{q}-|v|^q)\|_{L^2}
&\les
 &
(\|\chi_2 u\|_{L^{2q}_rL^\infty_\omega}+\|\chi_2 v\|_{L^{2q}_rL^\infty_\omega})^{q-1}
\|\chi_2 (u-v)\|_{L^{2q}_r L^2_\omega}\nonumber\\
&\les&
\langle t\rangle^{\al_4}(Z_u^2(t)+Z_v^2(t))^{q-1}Z_{u-v}(t),
\end{eqnarray}
which, together with \eqref{eq-claim4'-prf}, gives \eqref{eq-claim4'}.


\subsection{Proposition \ref{thm-data} }
In this subsection, we give the tedious proof for Proposition \ref{thm-data}.


By simple calculation, we can easily see that 
\begin{align}
 \Gamma^{\alpha}Pu(0)
 & = \sum_{|b| \leq 2}\sum_{|a| \leq |b|} C^{\alpha}_{ab}x^a \partial^{b} Pu(0), \\
 \partial_t\Gamma^{\alpha}Pu(0)
 & = \sum_{1\leq |b|\leq 2}\sum_{|a| \leq |b|-1} \tilde{C}^{\alpha}_{ab}x^a 
 \partial^{b} Pu(0) + \sum_{|a| \leq |b| = 2} \tilde{C}^{\alpha}_{ab}x^a 
 \partial^{b} \partial_t Pu(0)
\end{align}
for $|\alpha| \leq 2$, where $x^a = x_1^{a_1}\cdots x_n^{a_n}$ and 
$\partial^b = \partial_t^{b_0}\cdots \partial_n^{b_n}$. 
Thus we have
\begin{eqnarray}
&&\label{y1}
\sum_{|\alpha|\leq 2}
\|\pa \Gamma^\alpha P u(0)\|_{L^2\cap {\dot H}^{s_d-1}}
\\
&\les &\sum_{\substack{|a| \leq |b| \leq 2}} \left\|x^a 
 \partial^{b} Pu(0)\right\|_{\dot{H}^1\cap \dot{H}^{s_d}}
+\sum_{\substack{|a| \leq |b|-1\le 1 }} \left\|x^a 
 \partial^{b} Pu(0)\right\|_{L^2\cap\dot{H}^{s_d-1}}\nonumber\\
&&+\sum_{|a| \leq |b| = 2} \left\|x^a 
 \partial^{b} \partial_t Pu(0)\right\|_{L^2\cap\dot{H}^{s_d-1}},\nonumber\\
&\les& \varepsilon\Lambda +
\sum_{|a| \leq 1} \left\|x^a 
 \partial_t^2 Pu(0)\right\|_{L^2\cap\dot{H}^{s_d-1}}
+\sum_{|a| \leq 2} \left\|x^a 
 \pa \partial_t^2 Pu(0)\right\|_{L^2\cap\dot{H}^{s_d-1}}\nonumber\\
&\les&\varepsilon\Lambda +\sum_{|a| \leq 1} \left\|x^a 
 \Box Pu(0)\right\|_{L^2\cap\dot{H}^{s_d-1}}
+\sum_{|a| \leq 2} \left\|x^a 
 \partial \Box Pu(0)\right\|_{L^2\cap\dot{H}^{s_d-1}}.\nonumber
\end{eqnarray}
Therefore, the proof of \eqref{eq-data} is reduced to the estimate
\begin{equation}
\label{y2}
 \sum_{\substack{|a| \leq |b|+1 \leq 
  2}}\left\|x^a\partial^b\Box Pu(0)\right\|_{L^2 \cap \dot H^{s_d-1}}
 \les \varepsilon \left(\Lambda^p + \Lambda^q + \Lambda^{2p-1}+\Lambda^{p+q-1}\right).
\end{equation}

According to the definition of 
$Pu$, we have
\begin{eqnarray}
&&\sum_{|a|\le |b|+1 \leq 2}
\|x^a \partial^b\Box Pu(0)\|_{L^2  \cap \dot H^{s_d-1}} 
\label{y3}\\
 &\les & 
 \sum_{|a|\le 1 }
\|x^a (|\varepsilon f|^q+|\varepsilon g|^p)\|_{L^2  \cap \dot H^{s_d-1}} 
+  \sum_{|a|\le 2 }\|x^a |\varepsilon f|^{q-1}|\varepsilon (\nabla f,g)|\|_{L^2  \cap \dot H^{s_d-1}} \nonumber\\
 &
&+ 
\|x^{\le 2} |\varepsilon g|^{p-1}|\varepsilon \nabla (\nabla f, g)|\|_{L^2  \cap \dot H^{s_d-1}} 
+  \|x^{\le 2} |\varepsilon g|^{p-1}(|\varepsilon f|^q+|\varepsilon g|^p)\|_{L^2  \cap \dot H^{s_d-1}} \ .\nonumber
\end{eqnarray}
Recalling \eqref{eq-data1}, we observe that
\beeq\label{prop-data0}
\|x^k \nabla^l (\nabla f,g)\|_{\dot H^{2-l}\cap \dot H^{s_d+k-l-1}}\les \La, 0\le k\le 2, 0\le l\le 1\ ,\eneq
\beeq\label{prop-data1}\|x^k f\|_{\dot H^3\cap \dot H^{s_d+k-1}}+\|f\|_{\dot H^3\cap \dot H^{s_d}}
+ \|x^2 (\nabla f,g)\|_{\dot H^{1}\cap \dot H^{s_d}}\les \La, 1\le k\le 2\ .\eneq
In particular,
by Sobolev embedding, we have (see \eqref{eq-def-p_1}-\eqref{eq-rel})
$$
x^{\le 1} f\in \dot H^3\cap \dot H^{s_d}\subset L^\infty\cap L^{q_1},\  
g\in \dot H^2\cap \dot H^{s_d-1}\subset L^\infty\cap L^{2}\ .$$
$$
x^{\le 2}\nabla^{\le 1}(\nabla f, g)\in \dot H^1\cap \dot H^{s_d},\ 
x^{\le 2} g\in \dot H^1\cap \dot H^{s_d}\ ,$$
and for $n\ge 2$, $q\ge 2$,
$$\|uv\|_{L^2\cap \dot H^{s_d-1}}
\les
\|u\|_{L^n\cap L^{n+}}
\|v\|_{\dot H^1\cap \dot H^{s_d}}\ .$$
On the basis of these information, we are ready to prove \eqref{y2}. Actually, the first term on the right of \eqref{y3}, $\|x^{\le 1} (|\varepsilon f|^q+|\varepsilon g|^p)\|_{L^2  \cap \dot H^{s_d-1}}$, can be controlled by
$$
 \ep^q \|x^{\le 1} f\|_{\dot H^1\cap \dot H^{s_d}}
\|f\|^{q-1}_{L^{n(q-1)}\cap L^\infty}
+\ep^p \|x^{\le 1} g\|_{\dot H^1\cap \dot H^{s_d}}
\|g\|^{p-1}_{L^{n(p-1)}\cap L^\infty}\les \ep^q\La^q+\ep^p\La^p\ ,
$$
where we have used the fact that $n(p-1)\ge 2$ and $n(q-1)\ge q_1$ for
$n, p, q\ge 2$.
For the second and third terms, they are bounded by
$$
 \ep^q \|x^{\le 2} (\nabla f,g)\|_{\dot H^1\cap \dot H^{s_d}}
\|f\|^{q-1}_{L^{n(q-1)}\cap L^\infty}
+\ep^p \|x^{\le 2}\nabla(\nabla f,g)\|_{\dot H^1\cap \dot H^{s_d}}
\|g\|^{p-1}_{L^{n(p-1)}\cap L^\infty}\les \ep^q\La^q+\ep^p\La^p\ .
$$
For the last term,
\begin{eqnarray*}
&& \|x^{\le 2} |\varepsilon g|^{p-1}(|\varepsilon f|^q+|\varepsilon g|^p)\|_{L^2  \cap \dot H^{s_d-1}}\\
& \les&
 \ep^{p+q-1}\|x^{\le 2}g\|_{\dot H^1\cap \dot H^{s_d}} 
 \| g\|_{L^\infty}^{p-2}\|f\|_{L^{nq}\cap L^\infty}^q
 +
 \ep^{2p-1}
 \|x^{\le 2}g\|_{\dot H^1\cap \dot H^{s_d}} 
 \| g\|_{L^{n(2p-2)}\cap L^\infty}^{2p-2}\\
 & \les&
(\ep \La)^{p+q-1}+ (\ep \La)^{2p-1}\ .
\end{eqnarray*}
Hence \eqref{y2} is proved and the proof of \eqref{eq-data} is completed.

\section{Global existence}\label{sec-4}
In this section, using Propositions \ref{thm-unifbd} and
\ref{thm-conv}, we give the proof of global existence for 
$p>p_c$, $q>q_c$ and $((n-1)p-2)(q-1)\ge 4$, i.e., Theorem \ref{thm-1}.

With more notations, we could state a more precise version of the theorem.
\begin{theorem}\label{thm-1-2}
Let $n=2, 3$,
$s_d:=1/2-1/q$,
$q>q_c$, $p>p_c$ and $(q-1)((n-1)p-2)\ge 4$.
Suppose that $f\in \dot H^1\cap \dot H^{s_d}$ and $g\in L^2\cap \dot H^{s_d-1}$ with $\La<\infty$.
Then, there exists an $\varepsilon_0>0$ depending on 
$n$, $p$, $q$, and $\Lambda$ such that 
the Cauchy problem \eqref{nlw}-\eqref{eq-data0} admits a unique global solution,
provided that $\varepsilon\in [0,\varepsilon_0)$. Moreover, there exists a constant $C>0$ such that the solution satisfies the following estimates, 
$$X_{u}^2(t)\le CM\ep A(t), Y_{u}^2(t)\le CM\ep,$$
where
\beeq\label{grow}A(t)=\left\{
\begin{array}{ll}
1    ,  &   q>\max(q_c, \frac{2}{(n-1)(p-p_c)}), \\
\ln(2+t),&q=\frac{2}{(n-1)(p-p_c)}, p<q_c,  \\
   \<t\>^{\frac 1q-\frac{n-1}2 (p-p_c)},   &  1+\frac{4}{(n-1)p-2} \le q<\frac{2}{(n-1)(p-p_c)}, p<q_c.
\end{array}
\right.\eneq
See Figure \ref{3D-bu3} for an illustration of the region division.
\end{theorem}

As we see in the statement, we will give the proof of global existence for the following three cases:
\begin{enumerate}
  \item $p>p_c$, $q>\max(q_c, \frac{2}{(n-1)(p-p_c)})$,
  \item $p\in (p_c, q_c)$, $q=\frac{2}{(n-1)(p-p_c)}$,
  \item $p\in (p_c, q_c)$, $q\in [1+\frac{4}{(n-1)p-2} ,\frac{2}{(n-1)(p-p_c)})$,
\end{enumerate} 
and we will try to solve \eqref{nlw} in the ball  $B_{R}:=\{u\in S_\infty: 
X_{u}^2(t)\le R A(t), Y_{u}^2(t)\le R, \forall t\ge 0
\}$, equipped with a weaker topology defined by $\|u\|=\|A(t)^{-1} X_u(t)+Y_u(t)\|_{L^\infty}$.

 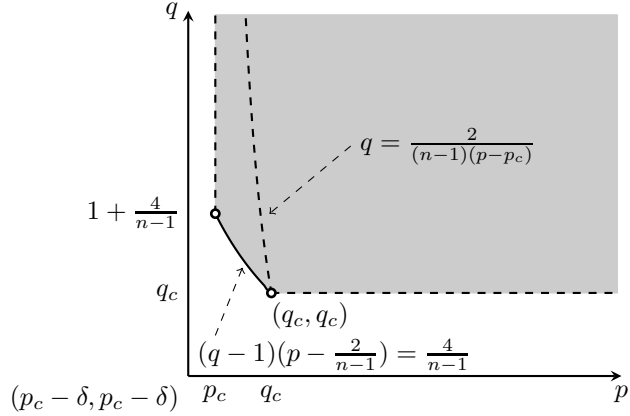
\begin{figure}
\centering
\begin{tikzpicture} [scale=1.8]
\filldraw[black!20!white] (4.97,4.47)--(2,4.47)--(2,3) to[out=297,in=134] (2.414,2.414)--(4.97,2.414)--(4.97,4.47);
\draw[thick, domain=2:2.414] plot (\x, {{1+2/((\x)-1)}}) ;
\draw[thick, dashed, domain=2.225:2.414] plot (\x,  {{1/((\x)-2)}});
\node[below] at (2.7,2.414) {$(q_c,q_c)$}; 
\node[left] at (1.8,3) {$1+\frac{4}{n-1}$};

\node[below left] at (1.8,1.8) {$(p_c-\de,p_c-\de)$};
\node[left] at (1.8,2.414) {$q_c$};
\node[below] at (2.414,1.8) {$q_c$};
\node[below] at (2,1.8) {$p_c$};
	\draw[thick, dashed] (2,3)--(2,4.47) (2.414,2.414)--(4.97,2.414);
\draw[fill=white,line width=1pt] (2.414,2.414) circle[radius=0.3mm];
\draw[fill=white,line width=1pt] (2,3) circle[radius=0.3mm];
	\draw[thick,-stealth] (1.8,1.8)--(1.8,4.5) node[left]{$q$};
	\draw[thick,-stealth] (1.8,1.8)--(5,1.8) node[below]{$p$};
%
\draw[->, dashed] (2,2.1)--(2.2,2.6);
\node[right] at (1.8,1.95) {$(q-1)(p-\frac{2}{n-1})=\frac{4}{n-1}$};

\draw[->, dashed] (3,3.5)--(2.4,3);
\node[right] at (3,3.5) {$q=\frac{2}{(n-1)(p-p_c)}$};

\end{tikzpicture}
\caption{Region division in the region of global existence}
\label{3D-bu3}
\end{figure}

\subsection{Case 1}
Let us begin with the easier case:  $p>p_c$, $q>\max(q_c, \frac{2}{(n-1)(p-p_c)})$. 
It turns out that there exists an $\ep_0>0$, such that, for any $\ep\le \ep_0$
\beeq\label{eq-glob-1}
u\in B_{2 C_1M\ep} \Rightarrow
Pu\in B_{2 C_1M\ep}\ ,
\eneq
\beeq\label{eq-glob-2}
u,v\in B_{2 C_1M\ep}
\Rightarrow
\|Pu-Pv\|\le \frac{1}{2} \|u-v\|\ .
\eneq
On the basis of these two estimates, it is a standard process to get the global existence of a unique solution in $B_{2 C_1M\ep_0}$.

First, let us give the proof of \eqref{eq-glob-1}.
Here, observe by direct calculation that
$$\al_1<-1 \Leftrightarrow q>q_c,\ 
\al_3<-1 \Leftrightarrow p>p_c,\ 
\al_4<-1 \Leftrightarrow q>p_c,
$$
$$\al_2<-1 \Leftrightarrow q>\frac{2}{(n-1)(p-p_c)},\ 
\al_5<-1 \Leftrightarrow q>1+\frac 1{n-1}.
$$
Then we see that $\al_j<-1$ for all $1\le j\le 5$ for this range of $(p,q)$, and so, recalling \eqref{eq-key-1'}-\eqref{eq-key-2'}, we have
$$\max(X_{Pu}^2(T),Y_{Pu}^2(T))
\le  C_1 M\ep+C \ep^p+C\ep^q\le  2 C_1 M\ep\ ,$$
for any $T\ge 0$ and $\ep\le \ep_0$, with $\ep_0\le 1$ satisfying
$C \ep_0^{p-1}+C\ep_0^{q-1}\le  C_1 M$.

Similarly, for \eqref{eq-glob-2}, 
recalling \eqref{eq-key-3}-\eqref{eq-key-4},
we have 
$$
u, v\in B_{2 C_1M\ep}\Rightarrow
\|Pu-Pv\|
\le C (\ep^{q-1}+\ep^{p-1})\|u-v\|\le \frac 12 \|u-v\|,$$
with sufficiently small $\ep_0\ll 1$.

\subsection{Case 2}
For the ``critical" case:  $p\in (p_c, q_c)$, $q=\frac{2}{(n-1)(p-p_c)}$, 
comparing with the proof for case 1, the only difference is that we have $\al_2=-1$.
On the basis of this observation, we have 
$Y_{Pu}^2(T)\le 2C_1 M\ep$, and
$$X_{Pu}^2(T)
\le  C_1 M\ep+C \ep^p \ln (2+T)+C\ep^q\le  2 C_1 M\ep \ln (2+T)\ ,$$
$$X_{Pu-Pv}(T)\le C \ep^{p-1} \|u-v\|\ln (2+T) +C\ep^{q-1} \|u-v\|\le  \frac{1}{2} \|u-v\|\ln (2+T) \ ,$$
$$Y_{Pu-Pv}(T)\le C \ep^{p-1} \|u-v\|+C\ep^{q-1} \|u-v\|\le  \frac{1}{2} \|u-v\|\ ,$$
for any $T\ge 0$ and $\ep\le \ep_0$, with $\ep_0\ll 1$.

\subsection{Remaining cases}
Inspired by the proof for the ``critical" case 2, for the remaining cases
$p\in (p_c, q_c)$, $q\in [1+\frac{4}{(n-1)p-2} ,\frac{2}{(n-1)(p-p_c)})$,
we use the ansatz $X_u^2(t)\les \ep \<t\>^\ga$, $Y_u^2(t)\les \ep$ to give the proof.

Let $\ga=\frac 1q+1-\frac{n-1}2 (p-1)=\frac 1q-\frac{n-1}2 (p-p_c)>0$, we observe by direct calculation that $\al_2= \ga-1$,
$$
\al_1+\ga q\le \ga-1\Leftrightarrow 
\ga\le \frac{n-1}2-\frac{q+1}{q(q-1)}
\Leftrightarrow 
q\ge 1+\frac{4}{(n-1)p-2} ,
$$
$$\ \al_3<-1 \Leftrightarrow p>p_c,$$
$$\al_4+\ga q \frac{q+1}{q+2}<-1 \Leftrightarrow 
\ga<\frac{n-1}2\frac{(q+2)(q-p_c)}{q(q+1)} 
 \Leftrightarrow 
\frac{n-1}2 (p-1)>\frac{2}q+\frac{n}{q(q+1)}+\frac{3-n}{2} ,
$$
$$
\al_5+\ga q <-1 
\Leftrightarrow 
\ga<\frac{n-1}2-\frac{n}{2q}
\Leftrightarrow 
\frac{n-1}2 (p-1)>\frac{n+2}{2q}+\frac{3-n}{2}\ .
$$
Notice that the last two inequality are true when
$q=1+\frac{4}{(n-1)p-2} $ and so for any $q\ge 1+\frac{4}{(n-1)p-2}$.
This proves
$Pu\in B_{2C_1 M\ep}$ for any $u\in B_{2C_1 M\ep}$,
by \eqref{eq-key-1'}-\eqref{eq-key-2'}. In the same vein, we have the convergence.

\section{Long time existence}\label{sec-5} 
In this section, using Propositions \ref{thm-unifbd} and
\ref{thm-conv}, assuming
$p, q\ge 2$ and $q>2/(n-1)$,
we give the proof of long time existence, for the cases
$p\le p_c$, $q\le q_c$ or $((n-1)p-2)(q-1)<4$,
 i.e., Theorem \ref{thm-2}.

With more notations, we state a more precise version of the existence theorem in the following
\begin{theorem}\label{thm-2-2}
Let $n=2, 3$, $q>2/(n-1)$ and $q, p\ge 2$. Assume also $q\le q_c$, $p\le p_c$ or $(q-1)((n-1)p-2)<4$.
Then for any
$f\in \dot H^1\cap \dot H^{s_d}$ and $g\in L^2\cap \dot H^{s_d-1}$ with $\La<\infty$,
there exists an $\varepsilon_0>0$ depending on 
$n$, $p$, $q$, and $\Lambda$ such that 
the Cauchy problem \eqref{nlw}-\eqref{eq-data0} admits a unique solution, for $t\in [0, T]$,
provided that $\varepsilon\in (0,\varepsilon_0)$, where
$$T= \left\{
\begin{array}{ll}
G_\ep (p) ,    &    2\le p\le p_c,\ q\ge 2p-1,\\
S_\ep (q)    , &    2\le q< q_c,\ 2/(n-1)<q\le p,\\
\exp(c\ep^{1-q}) ,    &    q=q_c\le p,\\
Z_\ep (p,q) , & (q-1)((n-1)p-2)<4,\ 2\le p\le q\le 2p-1
\end{array}\right.
$$
for some small constant $c>0$.
Moreover, there exists a constant $C>0$ such that the solution satisfies the following estimates
$$X_{u}^2(t)\le CM\ep \<t\>^\ga, Y_{u}^2(t)\le CM\ep, t\in [0, T],$$
where
\beeq\label{grow2}\ga=\left\{
\begin{array}{ll}
0    ,  &   2< q\le q_c, p\ge q, \\
\frac1q,& 2\le p\le p_c, q\ge 2p-1,  \\
\frac{q+1}{q-1}(\frac 1p-\frac 1q), 
  & 
  ((n-1)p-2)(q-1)<4, p\le q\le 2p-1.\end{array}
\right.\eneq
\end{theorem}

We will give the proof of long time existence for the following cases:
\begin{enumerate}
    \item $p\ge q$, $\ga=0$, $T_\ep=
    \left\{
    \begin{array}{ll}
    S_\ep(q),&   2\le q<q_c , q>2/(n-1),\\
     \exp(c \ep^{-(q-1)}) ,&    q=q_c,
\end{array}\right.$
\item $q\ge 2p-1$, $2\le p\le p_c$, $\ga=1/q$, $T_\ep=G_\ep(p)$,
  \item $q\in (p, 2p-1)$, $((n-1)p-2)(q-1)<4$,
$\ga=\frac{q+1}{q-1}(\frac 1p-\frac 1q)$,
   $T_\ep=Z_\ep(p,q)$.
\end{enumerate} 
As in 
Section \ref{sec-4}, we will try to solve \eqref{nlw} in the ball  $B_{R}:=\{u\in S_{T_\ep}, 
X_{u}^2(t)\le R \<t\>^\ga, Y_{u}^2(t)\le R, \forall t \in [0, T_\ep]
\}$, equipped with a weaker topology defined by $\|u\|=\|\<t\>^{-\ga} X_u(t)+Y_u(t)\|_{L^\infty_{T_\ep}}$.

Before proceeding, let us state a technical lemma, which proof is elementary and is left for the interested readers.
\begin{lemma}\label{lem-key}
Let $T_\ep=c\ep^{-A}$ with $A>0$, $\ga\ge 0$ and $s>1$.
If $A(1+\be-\ga)\le s-1$,
then for any $\de>0$, there exist $c>0$ and $\ep_0>0$ such that
$$\int_0^T\<t\>^\be \ep^s dt\le \de \ep \<T\>^\ga, \forall T\in [0,T_\ep], \ep\in (0,\ep_0]\ .$$
Similarly, let $\ln T_\ep=c\ep^{-A}$ with $A>0$, and $s>1$. Then 
for any $\de>0$, there exist $c>0$ and $\ep_0>0$ such that we have
$$\int_0^T\<t\>^\be \ep^s dt\le \de \ep , \forall T\in [0,T_\ep], \ep\in (0,\ep_0]\ ,$$
if $\be<-1$ or $A\le s-1$ when $\be=-1$,
and
$$\int_0^T\<t\>^\be \ep^s dt\le \de \ep \<T\>^\ga, \forall T\in [0,T_\ep], \ep\in (0,\ep_0]\ ,$$
if $1+\be\le \ga$ and $\ga>0$.
Moreover,
$$\int_0^T\<t\>^\be \ep^s dt\le \de \ep \<T\>^\ga, \forall T\in [0,\infty), \ep\in (0,\ep_0]\ ,$$
if $1+\be\le \ga$ and $\ga\ge 0$, except the case of $1+\be=\ga=0$.
For the critical case $\be= -1$, we have
$$\int_0^T\<t\>^\be \ep^s dt\le \de \ep \ln (2+T), \forall T\in [0,\infty), \ep\in (0,\ep_0]\ ,$$
\end{lemma}

\subsection{Case 1}
Let us begin with the easier case:  $p\ge q>2/(n-1)$, $2\le q\le q_c$. 
Let
 $T_\ep=c \ep^{-\frac{2q(q-1)}{2(q+1)-(n-1)q(q-1)}}$
 when $q< q_c$, and   $T_\ep=\exp(c \ep^{-(q-1)})$ when $q=q_c$,
 with the constant $c$ to be determined.

Since $p\ge q$, we have
$$\al_1=\frac 1q-\frac{n-1}2(q-1),
\al_2=\frac 1q-\frac{n-1}2(p-1)\le \al_1,$$
$$\al_3:=-\frac{n-1}2(p-1),\ 
\al_4:=-\frac{n-1}2(q-1)\in [\al_3,\al_1],$$
$$\al_5:=-q(\frac n2-s)+\frac n 2=-\frac 12-\frac{n-1}2(q-1)\le \al_4 .$$
We get from \eqref{eq-key-1'}-\eqref{eq-key-2'}  that, for any $u\in B_{2 C_1M\ep}$,
\beeq\label{eq-key1-1}
X_{Pu}^2(T)-
 C_1 M\ep\les \int_0^T (\<t\>^{\al_1}\ep^q + \<t\>^{\al_2}\ep^p )dt
\les \int_0^T \<t\>^{\al_1}\ep^q dt
\ ,\eneq
\beeq
Y_{Pu}^2(T)-
  C_1 M\ep\les \int_0^T (\<t\>^{\al_3}\ep^p +\<t\>^{\al_4}\ep^q+ \<t\>^{\al_5}
\ep^q )dt
\les \int_0^T \<t\>^{\al_1}\ep^q  dt
\ .
\label{eq-key1-2}
\eneq

Observe that for any $T\in [0, T_\ep]$,
$$q=q_c \Leftrightarrow \al_1=-1, A=q-1 \Rightarrow
\int_0^T \<t\>^{\al_1}\ep^q dt\les \ln (2+T) \ep^q\le c \ep\ll \ep\ ,$$
$$q<q_c \Rightarrow A=
\frac{2q(q-1)}{2(q+1)-(n-1)q(q-1)}>0, A(1+\al_1)= q-1 \Rightarrow
\int_0^T \<t\>^{\al_1}\ep^q dt\ll \ep\ .$$
Then we get 
$$
u\in B_{2 C_1M\ep}
\Rightarrow 
Pu\in  B_{2 C_1M\ep}, \forall T\in [0, T_\ep]\ .
$$
Similarly,
recalling \eqref{eq-key-3}-\eqref{eq-key-4},
we have 
$$
u, v\in B_{2 C_1M\ep}\Rightarrow
\|Pu-Pv\|
\ll \|u-v\|\ .$$
In summary, there exist $c_0>0$ and $\ep_0>0$ such that
$$
u, v\in B_{2 C_1M\ep}\Rightarrow
Pu\in B_{2 C_1M\ep}, \|P u-Pv\|\le \frac 12 \|u-v\|,$$
for any $c\in (0, c_0]$ and $\ep\in (0, \ep_0]$.

\subsection{Case 2}\label{sec-5.2}
For $q\ge 2p-1$ and $2\le p\le p_c$,
let $\ga=1/q$, and
$$T_\ep=\left\{
\begin{array}{ ll}
c \ep^{-\frac{2(p-1)}{2-(n-1)(p-1)}}:=c\ep^{-A_G} ,      &   p<p_c, \\
   \exp(c \ep^{-(p-1)}):=\exp(c\ep^{-A_G})   ,  &   p=p_c
\end{array}
\right.$$ with $c$ to be determined.

Observe that 
$$p=p_c, q\ge 2p-1 \Rightarrow 
\frac{n-1}2(q-1)\ge 2, \frac{n-1}2(p-1)=1
 \Rightarrow 
\al_1+\ga q+1\le \ga, \al_2+1=\ga,$$
$$\al_3=-1, \al_4+\frac{q+1}{q+2}q\ga<-1,
\al_5+q\ga< -1,$$
$$2\le p<p_c, q\ge 2p-1 \Rightarrow 
A_G(\al_1+\ga q+1- \ga)\le q-1,
A_G(\al_2+1- \ga)= p-1,$$
$$
A_G(\al_3+1)= p-1,
A_G(\al_4+\frac{q+1}{q+2}q\ga+1)\le q-1,
A_G(\al_5+q\ga+1)\le q-1.$$
By Lemma \ref{lem-key}, we have for any $T\in [0, T_\ep]$,
$$
\int_0^T \<t\>^{\al_1+\ga q}\ep^q dt
+\int_0^T \<t\>^{\al_2}\ep^p dt
\ll \ep \<T\>^\ga$$
$$\int_0^T \<t\>^{\al_3}\ep^p dt
+\int_0^T \<t\>^{\al_4+\frac{q+1}{q+2}q\ga}\ep^q
+\<t\>^{\al_5+q\ga}\ep^q
dt
\ll \ep \ .$$
Then, as a consequence of  \eqref{eq-key-1'}-\eqref{eq-key-4},
we conclude
that there exist $c_0>0$ and $\ep_0>0$ such that, 
$$
u, v\in B_{2 C_1M\ep}\Rightarrow
Pu\in B_{2 C_1M\ep}, \|P u-Pv\|\le \frac 12 \|u-v\|,$$
for any $c\in (0, c_0]$ and $\ep\in (0, \ep_0]$.

\subsection{Case 3}
In this case, we set
$\ga=\frac{q+1}{q-1}(\frac 1p-\frac 1q)=\frac{q+1}{q-1}\frac 1p+\frac 1q-\frac{2}{q-1}$,
and   $T_\ep=c \ep^{-\frac{2p(q-1)}{2(q+1)-(n-1)p(q-1)}}$ with $c$ to be determined.

As $((n-1)p-2)(q-1)<4$,
we have $A:=\frac{2p(q-1)}{2(q+1)-(n-1)p(q-1)}>0$, then $T_\ep=c \ep^{-A}$ and
$\frac{1}A=\frac{q+1}{p(q-1)}-\frac{n-1}2$.

Observe that 
$$
A(\al_1+\ga q+1- \ga)= q-1,
A(\al_2+1- \ga)= p-1\ ,$$
$$
q\le 2p-1\Leftrightarrow A(\al_3+1)\le p-1\ ,
$$
$$p\ge 2\frac{q+1}{q+3}\Leftrightarrow
A(\al_4+\frac{q+1}{q+2}q\ga+1)\le q-1,
A(\al_5+q\ga+1)\le q-1\ .$$
Then by Lemma \ref{lem-key}, we have for any $T\in [0, T_\ep]$,
$$
\int_0^T \<t\>^{\al_1+\ga q}\ep^q dt
+\int_0^T \<t\>^{\al_2}\ep^p dt
\ll \ep \<T\>^\ga\ ,$$
$$\int_0^T \<t\>^{\al_3}\ep^p dt
+\int_0^T (\<t\>^{\al_4+\frac{q+1}{q+2}q\ga}\ep^q
+\<t\>^{\al_5+q\ga}\ep^q)
dt
\ll \ep \ .$$
Thus, as a consequence of  \eqref{eq-key-1'}-\eqref{eq-key-4},
we conclude
that there exist $c_0>0$ and $\ep_0>0$ such that, 
$$
u, v\in B_{2 C_1M\ep}\Rightarrow
Pu\in B_{2 C_1M\ep}, \|P u-Pv\|\le \frac 12 \|u-v\|,$$
for any $c\in (0, c_0]$ and $\ep\in (0, \ep_0]$.


\bibliographystyle{amsplain}

%
\end{document}